\documentclass[11pt,reqno,sumlimits]{amsart}

\usepackage[margin=1in]{geometry}
\usepackage{amsmath, amsthm, amssymb}
\usepackage{times}
\usepackage[usenames,dvipsnames]{color}
\usepackage[colorlinks=true, pdfstartview=FitV, linkcolor=blue,
            citecolor=red, urlcolor=blue]{hyperref}
            
\usepackage{graphicx}

\date{\today}

\title[MG Survey]{The magneto-geostrophic equations: a survey}
\date{\today}
\thanks{This survey article was invited for the special volume of ``Proceedings of the St. Petersburg Mathematical Society'' dedicated to
O.A. Ladyzhenskaya. }

\author{Susan Friedlander}
\address{Department of Mathematics,
University of Southern California, 
Los Angeles, CA 90089} 
\email{\tt susanfri@usc.edu}

\author{Walter Rusin}
\address{Department of Mathematics,
University of Southern California, 
Los Angeles, CA 90089} 
\email{\tt wrusin@usc.edu}

\author{Vlad Vicol}
\address{Department of Mathematics,
Princeton University, 
Princeton, NJ 08544} 
\email{\tt vvicol@math.princeton.edu}

\usepackage{hyperref}

\theoremstyle{plain}
\newtheorem{theorem}{Theorem}[section]
\newtheorem{definition}[theorem]{Definition}
\newtheorem{lemma}[theorem]{Lemma}

\newtheorem{corollary}[theorem]{Corollary}

\theoremstyle{definition}
\newtheorem{remark}[theorem]{Remark}

\def\tilde{\widetilde}
\numberwithin{equation}{section}

\renewcommand\hat{\widehat}
\def\ZZ{{\mathbb Z}}

\def\RR{{\mathbb R}}
\def\TT{{\mathbb T}}
\def\PP{\mathcal P}
\def\QQ{\mathbb Q}

\def\Th{\Theta}

\def\supp{\mathop {\rm supp} \nolimits}

\def\eps{\varepsilon}

\begin{document}


\begin{abstract}
We discuss recent results obtained by the authors in~\cite{FriedlanderVicol11a,FriedlanderVicol11b,FriedlanderVicol12,FriedlanderRusinVicol12}, regarding the analysis of the magneto-geostrophic equation: a model proposed by Moffatt and Loper~\cite{MoffattLoper94,Moffatt08} to study the geodynamo and turbulence in the Earth's fluid core. We conclude this review by indicating some open problems around the MG equation, that remain to be addressed in the future.
\end{abstract}



\maketitle

\setcounter{tocdepth}{2} 
\tableofcontents

\section{Introduction}\label{sec:intro}

Physicists have long realized the importance of the Earth's magnetic field and that this field originates in the Earth's fluid core. The geodynamo is the process by which the rotating, convecting, electrically conducting molten iron in the Earth's fluid core maintains the geomagnetic field against ohmic decay. The convective processes in the core that produce the velocity fields required for dynamo action are a combination of thermal and compositional convection. The full dynamo problem requires the examination of the full three dimensional partial differential equations governing convective, incompressible magnetohydrodyamics (MHD). In the past decades computer models have been used to simulate the actual geodynamo. See, for example, the review articles by Roberts et al on modeling the geodynamo~\cite{BraginskyRoberts95,GlatzmaierRoberts97,RobertsGlatzmaier00,RobertsGlatzmaier01} and references therein. However, as Glatzmaier et al~\cite{GlatzmaierOgdenClune04} remark, no three dimensional dynamo model has been run at the spatial resolution required to simulate the broad spectrum of turbulence which surely exists in the Earth's fluid core. 

Current computers and numerical methods require the imposition of diffusivities that are several orders of magnitude larger than those which are realistic. It is therefore reasonable to attempt to gain some insight into the geodynamo by considering a reduction of the full MHD equations to a system that is more tractable, but one that retains many of the essential features of the problem. The magnetogeostrophic equation  proposed by Moffatt and Loper~\cite{MoffattLoper94, Moffatt08} is one such model. The physical postulates of this model are the following: slow cooling of the Earth leads to slow solidification of the liquid metal core onto the solid inner core and releases latent heat of solidification that drives compositional convection in the fluid core. The arguments for the appropriate ranges of the characteristic length, velocity, and perturbation density are based on these physical postulates.

We first present the full coupled three-dimensional MHD equations for the evolution of the velocity vector $u(x,t)$, the magnetic field vector $B(x,t)$ and the buoyancy field $\theta(x,t)$ in the Boussinesq approximation and written in the frame of reference rotating with angular velocity $\Omega$. The physical forces governing this system are Coriolis forces, Lorentz forces, and gravity. Following the notation of Moffatt and Loper~\cite{MoffattLoper94} we write the equations in terms of dimensionless variables. The orders of magnitude of the resulting non-dimensional parameters are motivated by the physical postulates of the Moffatt and Loper model:
\begin{align}
& N^2 \left[ R_o \left( \partial_t u + u\cdot \nabla u\right) +  {e}_3 \times u \right] = - \nabla P +  {e}_2 \cdot \nabla b + R_m b\cdot\nabla b + N^2 \theta  {e}_3 + \eps \Delta u \label{eq:1.1} \\
& R_m \left[ \partial_t b + u \cdot \nabla b - b\cdot \nabla u\right] =  {e}_2 \cdot \nabla u + \Delta b \label{eq:1.2} \\
& \partial_t \theta + u \cdot \nabla \theta = \eps_\kappa \Delta \theta + S \label{eq:1.3} \\
& \nabla \cdot u = 0, \nabla \cdot b = 0. \label{eq:1.4}
\end{align}
Here $S(x,t)$ is a given smooth function that represents the forcing of the system.

The mathematical statement of the geodynamo problem asks whether there are initial data for the MHD system for which the evolution of the perturbation of the magnetic field $b(x,t)$ grows for sufficiently long time. This can be interpreted as a question of the existence of instabilities of \eqref{eq:1.1}--\eqref{eq:1.4}.

The notation in \eqref{eq:1.1}--\eqref{eq:1.4} is the following. The Cartesian unit vectors are denoted by $ {e}_1, {e}_2$, and $ {e}_3$. For simplicity, we have assumed that the axis of rotation and the gravity $g$ are aligned in the direction of $ {e}_3$. We have assumed that the magnetic field $B(x,t)$ consists of the sum of an underlying constant field $B_0  {e}_2$ and a perturbation $b(x,t)$. Our choice of $B_0 e_2$ as the underlying magnetic field is consistent with the models where the magnetic field is believed to be predominantly toroidal (cf.~\cite{MoffattLoper94}). The dimensionless parameters are
\begin{align}
N^2 = \frac{2 \Omega \mu_0 \eta \rho}{B_0^{2}}, \quad 
R_o = \frac{V}{2 L \Omega}, \quad
R_m = \frac{VL}{\eta}, \quad
\eps = \frac{\nu \eta \mu_0 \rho}{B_0^2 L^2}, \quad
\eps_\kappa = \frac{\kappa}{LV} \label{eq:parameters}
\end{align}
where the typical velocity scale $V$ is chosen to be 
\begin{align}
V = \frac{\theta_0 g}{2\Omega}. \label{eq:V}
\end{align} 
Here $\theta_0$ is a typical amplitude of the buoyancy $\theta$ and $L$ is a characteristic length scale. $P$ is the sum of the fluid and magnetic pressures, $\nu$ is the kinematic viscosity, $\eta$ is the magnetic diffusivity, and $\kappa$ is the molecular diffusivity of the compositional variation that creates an ambient density $\rho$. The usual estimate for $\eta$ in the fluid core is $3 \mbox{m}^2\; \mbox{cm}^{-1}$. The values of $\nu$ and $\kappa$ are speculative, but likely to be extremely small. For a more detailed discussion of plausible ranges of the physical parameters that are appropriate for the geodynamo, we refer the reader to the book of Ghil and Childress~\cite{GhilChildress87}.

For the regions in the Earth's fluid core modeled by Moffatt and Loper, it is argued that the dimensionless parameters are the following orders of magnitude. The inverse Elsasser number $N^2$ will be of order unity, or somewhat less in the condition of the fluid core, while the Rossby number $R_o$ is of order $10^{-3}$. The value of $\eps$ is somewhat uncertain since the viscosity $\nu$ is a matter of guess work, but plausibly $\eps \le 10^{-8}$ in the core conditions. The inverse Peclet number $\eps_\kappa$ is also likely to be extremely small, perhaps of the order of $10^{-8}$. According to Moffatt and Loper~\cite{MoffattLoper94}, the magnetic Reynolds number $R_m$ is relatively small in the range of $L$ considered. With these estimates of the sizes of the parameters the dominant balance of the leading order terms in equations \eqref{eq:1.1} and \eqref{eq:1.2} give the following reduced system:
\begin{align}
& N^2  {e}_3 \times u = - \nabla P +  {e}_2 \cdot \nabla b + N^2 \theta  {e}_3 \label{eq:1.5}\\
& 0 =  {e}_2 \cdot \nabla u + \Delta b \label{eq:1.6}
\end{align}
along with the incompressibility condition
\begin{align}
\nabla \cdot u = 0, \nabla \cdot b = 0 \label{eq:1.7}.
\end{align}
The {\em linear} system of equations \eqref{eq:1.5}--\eqref{eq:1.7} determine the differential operators that relate the vector fields $u(x,t)$ and $b(x,t)$ with the scalar buoyancy $\theta(x,t)$. These operators encode the vestiges of the physics in the problem, namely Coriolis force, Lorentz force, and gravity. Vector manipulations of \eqref{eq:1.5}--\eqref{eq:1.7} give the expression 
\begin{align}
\left[ N^4 ( {e}_3 \cdot \nabla)^2 \Delta + ( {e}_2 \cdot \nabla)^4 \right] u = N^4 ( {e}_3 \cdot \nabla) \Delta ( {e}_3 \times \nabla \theta) + N^2 ( {e}_2 \cdot \nabla)^2 \nabla \times ( {e}_3 \times \nabla \theta) \label{eq:1.8}.
\end{align}
The sole remaining nonlinearity in the system comes from the coupling of \eqref{eq:1.8} with the equation \eqref{eq:1.4} for the  time evolution of the $\theta(x,t)$. We call this nonlinear active scalar equation the {\em magnetogeostrophic} (MG) equation. It is the subject of this article.

We study the properties and behavior of the active scalar equation in $\TT^3$ 
\begin{align}
&\partial_t \theta + u\cdot \nabla \theta = \eps_\kappa \Delta \theta + S \label{eq:MG:1} \\
& u = M[\theta] \label{eq:MG:2}
\end{align}
via an examination of the Fourier multiplier operator $M$ obtained from \eqref{eq:1.8}, which relates $u$ and $\theta$. More precisely, 
\[ u_j = M_j[\theta] = ( \hat{M}_j \hat{\theta})^{\vee}\]
for $j \in \{1,2,3\}$, and the explicit expressions for the components of $\hat{M}$ as functions of the Fourier variable $k = (k_1,k_2,k_3) \in \ZZ^3_\ast:= \ZZ^3 \setminus \{ k_3=0\}$ are
\begin{align}
\hat{M}_1(k) &= \left( N^4 k_2 k_3 |k|^2 - N^2 k_1 k_2^2 k_3\right) D^{-1} \label{eq:M1}\\
\hat{M}_2(k) &= \left( - N^4 k_1 k_3 |k|^2 - N^2   k_2^3 k_3\right) D^{-1} \label{eq:M2}\\
\hat{M}_3(k) &= \left( N^2 k_1^2 k_2^2 + N^2  k_2^4 \right) D^{-1} \label{eq:M3}
\end{align}
where
\begin{align}
D&= N^4 |k|^2 k_3^2 + k_2^4 \label{eq:D}.
\end{align}
On $\{k_3=0\}$ we define $\hat{M}_j(k) = 0$, for all $j\in\{1,2,3\}$, since for self-consistency of the model we assume that $\theta$ and $u$ have zero vertical mean, i.e. $\int \theta dx_3 = \int u dx_3 = 0$. Several important properties of the $\hat{M}_j$'s are immediately obvious. Firstly, the functions are strongly anisotropic with respect to the dependence on the integers $k_1, k_2$, and $k_3$. This is a consequence of the interplay of the three above mentioned physical forces. Secondly, the function $\hat{M}_3$ is an even function of $k_1, k_2$, and $k_3$. Thirdly, we have $\hat{M}_j(k) k_j = 0$ for any $k \in \ZZ^3$ (we use the standard summation convention on repeated indices), which is built into the system since the velocity field is divergence free
\[ \nabla \cdot u = 0.\]
Lastly, although the symbols $\hat{M}_j$ are bounded (by a constant that depends only on $N$) in regions of Fourier space where $k_2 \leq \max\{ k_2,k_3\}$, this is not the case in the ``curved'' regions where $k_3 = {\mathcal O}(1)$, and $k_2 = {\mathcal O}( |k_1|^r)$ with $0 < r \leq 1/2$.
In such regions the symbols are unbounded as $|k_1| \to \infty$ with 
\begin{align}\label{eq:M:max:growth}
|\hat{M}(k)| \leq C |k|,
\end{align}
for some suitable chosen constant $C$ that {\em is independent of $N$}, and this bound is sharp. Hence the MG equation falls in the class of abstract active scalar equations (see \eqref{eq:ABS:1}--\eqref{eq:ABS:2} below) where the relation between the divergence-free velocity field $u$ and the scalar $\theta$ is given as a singular operator of order $1$. These properties of $\hat{M}$ make the MG equation interesting and challenging mathematically, as well as having a clear physical basis in its derivation from the MHD equations.

\begin{figure}[!htb]
      \includegraphics[width=0.3\textwidth]{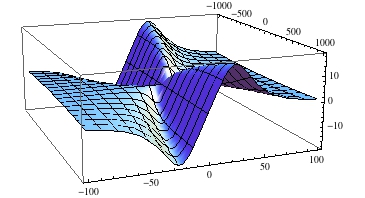} 
      \includegraphics[width=0.3\textwidth]{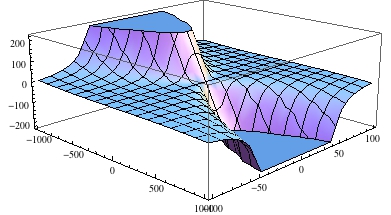} 
      \includegraphics[width=0.3\textwidth]{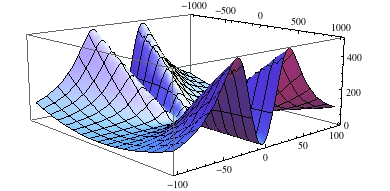}
      \caption{From left to right: plot of $\hat{M}_1(k_1,k_2,1)$, $\hat{M}_2(k_1,k_2,1)$, and $\hat{M}_3(k_1,k_2,1)$, with $N=1$.  The range of the frequencies is $k_1 \in [-10^3,10^3]$ and $k_2 \in [-10^2,10^2]$.}
\end{figure}

In Section~\ref{sec:diffusive} we state results proved in \cite{FriedlanderVicol11a,FriedlanderVicol11b,FriedlanderVicol12} for the diffusive MG equation, i.e., the case when $\eps_\kappa>0$. We recall that in the fluid core $\eps_\kappa$ is positive but probably extremely small. In this case the MG equation is global well posed. To prove this fact we use the control on $\Vert \theta(t) \Vert_{L^\infty}$, a quantity that is non-increasing along solutions, and which gives a bound on the drift velocity $u$ in the space $L^\infty_t BMO^{-1}_x$. The latter information on the drift is ``critical'' for the natural scaling of the {\em linear} drift-diffusion equation
\begin{align}
\partial_t \theta + v \cdot \nabla \theta = \Delta \theta \label{eq:linear}.
\end{align}
Indeed, for $v \in L^\infty_t BMO^{-1}_x \cap L^2_t L^2_x$  divergence-free, we proved in \cite{FriedlanderVicol11a} that
the solution $\theta$ of \eqref{eq:linear} becomes H\"older continuous in finite time, by appealing to the classical ideas of De Giorgi~\cite{DeGiorgi57}, in the spirit of~\cite{CaffarelliVasseur10}. This result was also established independently in~\cite{SereginSilvestreSverakZlatos12}. See also~\cite{Nash58,Osada87,Zhang04},  the more recent~\cite{NazarovUraltseva12,SilvestreVicol12,Filonov12,IgnatovaKukavicaRyzhik12}, and references therein, for regularity issues of linear parabolic equations with singular coefficients and divergence-free drift.
The H\"older continuity of solutions to the linear equation~\eqref{eq:linear} may then be bootstrapped to show that the solutions of the dissipative MG equations become smooth in arbitrarily short finite time~\cite{FriedlanderVicol11b}.

Having established the global well-posedness of \eqref{eq:MG:1}--\eqref{eq:MG:2} for any $\eps_\kappa>0$, we can proceed to examine the stability of equilibria. We show that the existence of an unstable eigenvalue for the linearized equations implies nonlinear (Lyapunov) instability of the nonlinear MG equation (see also~\cite{FriedlanderStraussVishik97,BardosGuoStrauss02,VishikFriedlander03,FriedlanderPavlovicShvydkoy06,FriedlanderPavlovicVicol09}). We then exhibit the existence of an unstable eigenvalue in a particular example, using the method of continued fractions cf.~\cite{MeshalkinSinai61,FriedlanderStraussVishik97}. The magnitude of this eigenvalue is ${\mathcal O}(\eps_\kappa^{-1})$ which, in the context of the fluid core, is extremely large. Thus the MG equation can be used to model very rapid exponential growth in time of $\theta$, and hence, through \eqref{eq:MG:1}--\eqref{eq:MG:2}, of the perturbation velocity and magnetic fields, i.e. a dynamo type instability. 
We emphasize that this ``MG-dynamo'' instability requires the presence of all three forces -- Coriolis, Lorentz, and gravity -- because the dimensionless parameter $N^2$ (defined above in \eqref{eq:parameters}) must be strictly in the range $(0,\infty)$.

In Section~\ref{sec:inviscid} we examine the non diffusive MG equation and set $\eps_\kappa=0$ (cf.~\cite{FriedlanderVicol11b}). The situation is then dramatically different from the case $\eps_\kappa>0$. Without the Laplacian to control the unbounded operator $M$, the MG equation is {\em ill-posed} in the sense of Hadamard in Sobolev spaces. Such a result utilizes in an essential way both the unboundedness of $\hat{M}_3$ and the fact that it is an even function of $k$. The proof again uses the existence of an unstable eigenvalue of the linearized equation, however when $\eps_\kappa=0$ it is possible to construct such an eigenvalue that is arbitrarily large and it then follows that there is no Lipschitz solution map at time $t=0$ (see also~\cite{Tao06,Renardy09,GerardVaretDormy10,GerardVaretNguyen10,GuoNguyen10}).

In Section~\ref{sec:frac} we describe results for the fractionally diffusive MG equation (cf.~\cite{FriedlanderRusinVicol12}), i.e. replace the Laplacian by $- (-\Delta)^\gamma$, for $\gamma \in (0,1)$. This situation is non-physical but mathematically interesting in the context of the tools required to treat the fractional Laplacian. We now have a dichotomy across the value $\gamma=1/2$. In the range $1/2 < \gamma < 1$ the equation is {\em locally} well-posed, while it is Hadamard {\em Lipschitz ill-posed} for $0<\gamma<1/2$. At the critical value $\gamma = 1/2$ the problem is globally well-posed for suitably small initial data (and source term), but is ill-posed for sufficiently large initial data. In terms of the non-dimensional parameters introduced above in \eqref{eq:parameters}, this dichotomy can be seen as global well-posedness for $\eps_\kappa \gg 1$, and Lipschitz ill-posedness for $\eps_\kappa \ll 1$. A further feature of interest is that the anisotropy of the symbol can be explored to obtain an improvement in the regularity of the solutions when the initial data and source are {\em supported on a plane in Fourier space}. For such well-prepared data the local existence and uniqueness of solutions can be obtained for all values $\gamma \in (0,1)$, and the global existence holds for all initial data when $\gamma \in (1/2,1)$.

We note that some of the theorems discussed in this article are valid for certain classes of partial differential equations cf.~\eqref{eq:ABS:1}--\eqref{eq:ABS:2} below, and in any spatial dimension $d\geq 2$. For such theorems, the three dimensional MG equation is a particular example of an active scalar equation for which these results hold.

We conclude the manuscript with Section~\ref{sec:concl}, in which we present a few remaining open problems regarding the mathematical analysis of the MG equations.

\section{Diffusive MG equation}
\label{sec:diffusive}
In this section we discuss the global regularity of solutions to \eqref{eq:MG:1}--\eqref{eq:MG:2} and the Lyapunov stability of steady states, in the case that $\eps_\kappa>0$.

\subsection{Global regularity of singular drift-diffusion equations}
\label{sec:global}
In order to show that there is no finite time loss of regularity of solutions to \eqref{eq:MG:1}--\eqref{eq:MG:2} with $\eps_\kappa>0$ let us first perform a rough derivative count. The dissipation is smoothing by two derivatives, while the nonlinearity is also losing two derivatives: the first coming from $\nabla \theta$, while the second one derivative loss is coming from the constitutive law $u = M[\theta]$, in view of \eqref{eq:M:max:growth}. In this sense we may call the problem ``critical'' (we make this more precise below).

Since the operator $M$ relating $u$ and $\theta$ is an order one 
pseudo-differenential operator, one may view the MG equation as an example of a more general class of of active scalar equations:
 \begin{align}
  &\partial_t \theta + (u\cdot \nabla) \theta = \Delta \theta \label{eq:ABS:1}\\
  &u_j = \partial_i T_{ij} \theta, \quad  \partial_j u_j = 0\label{eq:ABS:2}
\end{align}
with $T_{ij}$ a matrix of zero-order pseudo-differential operators, and we use the summation convention on repeated indices. 
We have dropped the smooth source term $S$ from \eqref{eq:MG:1} since it is not an obstacle to the well-posedness theory, have re-scaled the time variable so that the diffusivity coefficient is $1$ instead of $\eps_\kappa$ (this is not an issue since we will construct global in time solutions), and the $T_{ij}$ corresponding to \eqref{eq:MG:2} are given by $T_{ij} = - \partial_i (-\Delta)^{-1} M_j$. We call this class of active scalar equation ``singular'' since the drift velocity is given in terms of the advected scalar by a constitutive law which is losing derivatives. We refer the interested reader in to the work of Chae, Constantin, Cordoba, Gancedo, and Wu~\cite{ChaeConstantinCordobaGancedoWu12} and also~\cite{FriedlanderGancedoSunVicol12, MiaoXue12} for further discussions about the well-posedness theory of singular active scalar equations.

To make the idea of criticality more precise, we note that the active scalar equation \eqref{eq:ABS:1}--\eqref{eq:ABS:2} has a natural scaling invariance built into it: if $\theta(t,x)$ is a solution, then $\theta_\lambda(t,x) = \theta(\lambda^2 t, \lambda x)$ is {\em also} a solution, with corresponding drift given by $u_\lambda(t,x) = \lambda u(\lambda^2 t, \lambda x) = \partial_i T_{ij} \theta_\lambda(t,x)$. Hence $L^\infty(\RR^d)$ is the critical Lebesgue space with respect to the natural scaling of the equation. Additionally \eqref{eq:ABS:1} is a scalar equation, so that $L^p$ norms are non-increasing 
in time for $1 \leq p \leq \infty$, and hence the $L^\infty$ norm is the ``strongest'' a priori controlled norm along solutions. In this sense we say that \eqref{eq:ABS:1}--\eqref{eq:ABS:2} is a ``critical'' equation.

Another active scalar equation derived from the Navier-Stokes equations for a rapidly rotating fluid in a  thin shell is the surface quasi-geostrophic equation, introduced by Constantin, Majda, and Tabak~\cite{ConstantinMajdaTabak94} (see also~\cite{ConstantinCordobaWu01,CordobaCordoba04,Wu04,Kiselev10} and references therein), which in the critically dissipative case reads
\begin{align}
  &\partial_t \theta + (u\cdot\nabla)\theta = - (-\Delta)^{1/2} \theta, \label{eq:SQG:1}\\
  &u = \nabla^\perp {(-\Delta)^{-1/2}} \theta \equiv(R_2\theta,-R_1 \theta)\label{eq:SQG:2}
\end{align}
and $R_i$ represents the $i^{th}$ Riesz transform. There are strong analogies between the systems \eqref{eq:SQG:1}--\eqref{eq:SQG:2} and \eqref{eq:ABS:1}--\eqref{eq:ABS:2}. The velocity field is divergence free in both cases, and the $L^\infty$ norm is critical with respect to the natural scaling of the equation. The main difference is that in for SQG the dissipative operator is non-local, but the drift velocity is not more singular than the advected scalar.

Inspired by the De Giorgi-type proof of global regularity for the critically dissipative SQG equation of Cafarelli and Vasseur~\cite{CaffarelliVasseur10} (see also~\cite{KiselevNazarovVolberg07,KiselevNazarov09,ConstantinVicol12} for different proofs), to show that solutions of \eqref{eq:ABS:1}--\eqref{eq:ABS:2} do not blow up in finite time, we first look at the associated linear problem, where the drift velocity is {\em given} and obeys the bounds that are a priori available for weak solutions of the nonlinear problem: $u \in L^\infty(0,\infty;BMO^{-1}) \cap L^2(0,\infty; L^2)$. In~\cite{FriedlanderVicol11a} we use a De Giorgi iteration technique to prove that a weak solutions of the {\em linear} drift-diffusion equation
\begin{align}
&\partial_t \theta + (v\cdot \nabla) \theta = \Delta \theta \label{eq:LIN:1} \\
&\nabla \cdot v = 0, \quad v_j = \partial_i V_{ij}, \quad V_{ij} \in L^\infty(0,\infty;BMO) \cap L^2(0,\infty;H^1) \label{eq:LIN:2} 
\end{align}
instantly become H\"older continuous. More precisely, the result is:
\begin{theorem}[\cite{FriedlanderVicol11a}] \label{thm:DeGiorgi}
Let the initial data $\theta_0 \in L^2$ and $\theta \in L^\infty(0,\infty;L^2) \cap L^2(0,\infty;H^1)$ be a weak solution of \eqref{eq:LIN:1}, where the drift velocity $v$ is given by \eqref{eq:LIN:2} for a suitable matrix $V_{ij}$. Then, for any $t_0$, there exists $\alpha>0$ such that 
$\theta \in C^{\alpha/2}(t_0,\infty;C^\alpha)$.
\end{theorem}
The above result was obtained independently by Sergin, Silvestre, Sverak, and Zlatos in~\cite{SereginSilvestreSverakZlatos12}.
In view of the constitutive \eqref{eq:ABS:2} a direct corollary of Therem~\ref{thm:DeGiorgi} is that:
\begin{corollary}
Weak solutions $\theta$ of \eqref{eq:ABS:1}--\eqref{eq:ABS:2} are H\"older continuous in positive time, i.e. they obey $\theta\in L^\infty(t_0,\infty;C^\alpha)$ for any $t_0>0$, and a suitable $\alpha>0$.
\end{corollary}
In view of the above discussion about scaling, H\"older regularity is sub-critical information for \eqref{eq:ABS:1}-\eqref{eq:ABS:2} and one may bootstrap to prove that the solution is classical:
\begin{theorem}[\cite{FriedlanderVicol12}]\label{thm:higher}
Let $\theta \in L^\infty(t_0,\infty;C^\alpha)$ be a  weak solution of the initial value problem associated to \eqref{eq:ABS:1}--\eqref{eq:ABS:2}, with $\alpha \in (0,1)$. Then $\theta \in L^\infty(t_1,\infty;C^{1+\delta})$, for any $t_1 > t_0$, for some $\delta \in (0,1)$.
\end{theorem}
The proof of Theorem~\ref{thm:higher} relies on a suitable bootstrap procedure in Besov spaces, in the spirit of~\cite{ConstantinWu08}. Repeating the argument for the derivative of $\theta$ one may further obtain that the solution of \eqref{eq:ABS:1}--\eqref{eq:ABS:2} is in fact $C^\infty$ smooth in finite time. 

We conclude this subsection with a very brief discussion of the proof of Theorem~\ref{thm:DeGiorgi}. We refer the reader to~\cite{FriedlanderVicol11a} for all details.  Along the lines of Caffarelli and Vasseur~\cite[Theorem 3]{CaffarelliVasseur10} for the critical SQG equation, the first step consists of showing that a weak solution is bounded in the spacial variable (see~\cite[Appendix A]{FriedlanderVicol11a} for a construction of weak solutions which obey the energy inequality). This is achieved by proving that the function $(\theta-h)_+ = \max\{\theta-h,0\}$ has zero energy if $h$ is chosen large enough. The key ingredient is the inequality
\begin{align*}
  \partial_t (\theta - h)_+ - \Delta(\theta-h)_+ +(v\cdot\nabla)(\theta-h)_+ \leq 0
\end{align*}
which can be shown to hold for weak solutions. Using that $v$ is divergence free one obtains the energy inequality
\begin{align}\label{eq:L^2energy}
  \int_{\RR^d} | (\theta(t_2,\cdot) - h)_+|^2 dx + 2 \int_{t_1}^{t_2}\int{\RR^d} |\nabla (\theta - h)_+|^2 dx dt \leq \int_{\RR^d} | (\theta(t_1,\cdot) - h)_+|^2 dx,
\end{align}
for all $h>0$ and $0<t_1<t_2<\infty$. For $t_0 >0$, and $H>0$ to be chosen sufficiently large, letting $t_n = t_0 - t_0/2^n$, $h_n = H - H/2^n$, and
\begin{align*}
  c_n = \sup_{t\geq t_n} \int_{\RR^d} |(\theta(t,\cdot) - h_n)_+|^2 dx + 2 \int_{t_n}^{\infty} \int_{\RR^d} |\nabla (\theta-h_n)|^2 dx dt,
\end{align*}
it follows from \eqref{eq:L^2energy} that
\begin{align*}
  c_{n+1} \leq {\frac{C}{t_0 H^{4/d}}} 2^{n(1+4/d)} c_n^{1+2/d}
\end{align*}
for some constant $C>0$. Letting $H = H(\Vert \theta_0 \Vert_{L^2(\RR^d)},t_0^{d/4})$ be sufficiently large we thus obtain implies that $c_n \rightarrow 0$ exponentially as $n\rightarrow \infty$, and hence $ \theta(t_0,\cdot) \leq H$.

Having established the boundedness of the weak solutions, the next step is to prove H\"older regularity, which requires a more delicate argument. The main obstruction to applying the classical parabolic De Giorgi estimates via an $L^p$-based Caccioppoli inequality ($1<p<\infty$), is that $v_j = \partial_i V_{ij}$ is the derivative of a $BMO$ function. By essentially using the divergence-free nature of $v$, and appealing to the John-Nirenberg inequality, we may however prove a suitable local energy inequality: for any $0<r<R$ and $h\in \RR$, we have
\begin{align}\label{eq:1st:energyineq}
&\Vert (\theta-h)_+ \Vert_{L^{\infty}_t L^{2}_x(Q_r)}^2 + \Vert \nabla (\theta-h)_+ \Vert_{L^{2}_{t,x}(Q_r)}^2 \leq \frac{C\, R}{(R-r)^2} \Vert (\theta-h)_+ \Vert_{L^2_{t,x}(Q_R)}^{2-\frac{2}{d+2}} \Vert (\theta -h)_+ \Vert_{L^{\infty}_{t,x}(Q_R)}^{\frac{2}{d+2}},
\end{align}
where $C = C(d,\|V_{ij}\|_{L^\infty_t BMO_x})$ is a fixed positive constant. Here we have denoted by $Q_\rho$ the parabolic cylinder $[t_0 - \rho^2,t_0]\times B_\rho(x_0)$ for $\rho>0$ and $(t_0,x_0) \in (0,\infty)\times \RR^d$. Moreover, estimate \eqref{eq:1st:energyineq} also holds with $\theta$ replaced by $-\theta$. The bound \eqref{eq:1st:energyineq} may be in turn used, using a suitable iteration argument to obtain an improvement of oscillation result. Assume that {$h_0 \leq \sup_{Q_{r_0}}\theta$}, where $r_0 >0$ is arbitrary. Then we have
  \begin{align}\label{eq:sup:halfball}
  \sup_{Q_{r_0/2}}{\theta} \leq h_0 + C \left( \frac{|\{\theta > {h_0}\}\cap Q_{r_0}|^{1/(d+2)}}{r_0}\right)^{1/2} \left( \sup_{Q_{r_0}} {\theta} - h_0 \right)
  \end{align}
  for some positive constant $C = C(d,\|V_{ij}\|_{L^\infty_t BMO_x})$. We note that the prefactor of $\sup_{Q_{r_0}} {\theta} - h_0$ of the right side of \eqref{eq:sup:halfball} {\em does not scale with $r_0$}. Besides estimate \eqref{eq:sup:halfball}, the main ingredient in the proof of H\"older continuity is the following estimate on the growth of level sets of the solution. There exist $\kappa_0\in (0,1)$, $\delta_0 = \delta_0(\kappa_0, \|V_{ij}\|_{L^\infty_t BMO_x})$, and $n_0\geq 2$ such that for $t_1,R>0$ if we know that 
\begin{align*}
|\{ \theta(t_1,\cdot) \geq h\} \cap B_r | \leq \frac 12 |B_r|,
\end{align*}
then for all $t_2 \in [t_1,t_1 + \delta_0 r^2]$ we have
\begin{align}\label{eq:cond}
|\{ \theta(t_2,\cdot) \geq H\} \cap B_R | \leq \frac 78 |B_R|,
\end{align}
where $r=\kappa_0 R$, $M=\sup_{(t_1,t_1+\delta_0 R^2)\times B_R} {\theta}$, $m=\inf_{(t_1,t_1+\delta_0 R^2)\times B_R} {\theta}$, $h = (M+m)/2$, and $H = M-(M-m)/2^{n_0}$. To prove the above estimate one uses a suitable version of the local parabolic energy inequality. 

Bounds \eqref{eq:sup:halfball} and \eqref{eq:cond} are the main ingredients in the De Giorgi iteration procedure, which at the end of the day yields the improvement in oscillation 
\begin{align}
\left( \sup_{Q_1} \theta - \inf_{Q_1} \theta\right) \leq \gamma \left( \sup_{Q_2} \theta - \inf_{Q_2} \theta\right) \label{eq:osc}
\end{align}
for some $\gamma \in (0,1)$, where $Q_1 = [t_1,t_1+\delta_0 r^2] \times B_4$, $Q_2 = [t_1,t_1+\delta_0 R^2] \times B_R$, $t_1,R>0$ are arbitrary, $r = \kappa_0 R$, and $\delta_0, \kappa_0 \in (0,1)$ is as above.  Since $\gamma<1$, is then classical to show that \eqref{eq:osc} implies the solution $\theta$ is H\"older continuous, with exponent that may be computed explicitly.

\subsection{Lyapunov stability of steady states}
\label{sec:Lyapunov}

In this subsection we consider the stability of steady states to \eqref{eq:MG:1}--\eqref{eq:MG:2}. First we prove for steady states of the the more general class of equations \eqref{eq:MG:1}, with drift velocity given by \eqref{eq:ABS:2}, that {linear} implies {nonlinear} instability.  We then turn to the specific MG equations \eqref{eq:MG:1}--\eqref{eq:MG:2}, for which we construct an explicit unstable eigenvalue of the linear operator.

\subsubsection{Linear implies nonlinear instability}
Let $\Th_{0} \in C^{\infty}$ be a steady solution of \eqref{eq:MG:1}--\eqref{eq:ABS:2}, i.e. 
\begin{align} 
U_0 \cdot \nabla \Th_0 =  \eps_\kappa \Delta \Th_0 + S, \quad U_{0j} = \partial_i T_{ij} \Th_0, \quad \nabla \cdot U_0  = 0, \label{eq:steady:ABS}
\end{align}
where $S \in C^\infty$ is a time independent force, $\eps_\kappa>0$, and $T_{ij}$ are zero order pseudo differential operators.
Writing the solution of \eqref{eq:MG:1}--\eqref{eq:ABS:2} as $\theta + \Th_0$ we obtain that the perturbation temperature $\theta$ obeys
\begin{align}
&\partial_{t} \theta = L\theta + N\theta \label{eq:linearized},
\end{align}
where the {\em dissipative linear} operator $L$ is defined as
\begin{align}
L \theta = - U_{0} \cdot \nabla \theta - u \cdot \nabla \Th_{0} + \eps_\kappa \Delta \theta \label{eq:L}
\end{align}
with $u_{j} = \partial_{i} T_{ij}\theta$, and the {\em nonlinear} operator $N$ is given by
\begin{align}
N \theta  = - u \cdot \nabla \theta = - \nabla \cdot (u \theta).\label{eq:N}
\end{align}
Recall the following notion of Lyapunov stability:
\begin{definition}[~\cite{FriedlanderStraussVishik97}]\label{def:nonlinear:instability}
Let $(X,Z)$ be a pair of Banach spaces. A steady state $\Th_0$ is called $(X,Z)$ {nonlinearly stable} if for any $\rho>0$, there exists $\tilde{\rho}>0$ so that if
$\theta(\cdot,0) \in X$ and $\Vert{\theta(\cdot,0)}\Vert_{Z} < \tilde{\rho}$, then  
\begin{list}{\labelitemi}{\leftmargin=1em}
  \item[(i)] there exists a global in time solution $\theta$ to the initial value problem \eqref{eq:linearized}--\eqref{eq:N} with
  $\theta \in C([0,\infty);X)$;
  \item[(ii)] and we have the bound $\Vert{\theta(\cdot,t)}\Vert_{Z} < \rho$ for
  a.e.~$t\in [0,\infty)$.
\end{list}
An equilibrium $\Th_{0}$ that is not stable in the above sense is
called {Lyapunov unstable}.
\end{definition}
In the above definition, the Banach space $X$ is the space where a local existence theorem for the nonlinear equations is available (e.g.~$X = H^s$), while $Z$ is the space where the spectrum of the linear operator is analyzed, and where the instability is measured (e.g.~$Z=L^2$). The following result states that the linear Lyapunov instability of a steady state $\Th_0$ implies the nonlinear one.
\begin{theorem}[\cite{FriedlanderVicol11b}] \label{thm:linear:nonlinear}
 Let $\Th_0$ be a smooth, mean zero solution of \eqref{eq:steady:ABS}. If the associated linear operator
 $L$, as defined in \eqref{eq:L}, has spectrum in the unstable
 region, then the steady state is $(H^s,L^2)$ Lyapunov nonlinearly
 unstable, for arbitrary $s>d/2$.
\end{theorem}

The proof relies on a bootstrap argument that is by now well-known for fluids equations~\cite{FriedlanderStraussVishik97,BardosGuoStrauss02,VishikFriedlander03,FriedlanderPavlovicShvydkoy06,FriedlanderPavlovicVicol09}, and an a priori time independent bound on the $H^{2}$ norm of the solution, which we obtain from the De Giorgi $C^\alpha$ estimate. 

We only briefly sketch the argument given in full detail in \cite{FriedlanderVicol11b}, to which refer the interested reader. First, we fix a linearly unstable smooth eigenfunction $\phi$, with eigenvalue of {\em maximal} real part $\lambda$. In order to prove Theorem~\ref{thm:linear:nonlinear}, for $\eps \in (0,\eps_*]$ we consider solutions $\theta^\eps$ of
\begin{align} 
\partial_t \theta^\eps = L \theta^\eps + N(\theta^\eps), \quad \theta^\eps|_{t=0} = \eps \phi \label{eq:perturbation}
\end{align}
where $L,N,$ and $\phi$ are as above, and show that there exits a constant $C_* >0$ such that for any $\eps \in (0,\eps_*]$ we may find a time $T_\eps>0$ with 
\begin{align} 
\| \theta^\eps(T_\eps) \|_{L^2} \geq C_\ast. \label{eq:perturbation:lower}
\end{align}
By Definition~\ref{def:nonlinear:instability}, the bound \eqref{eq:perturbation:lower} shows that the trivial solution $\theta = 0$ of \eqref{eq:linearized}--\eqref{eq:N} is $(H^s,L^2)$ Lyapunov unstable: indeed, an ${\mathcal O}(\eps)$ perturbation leads in finite time to an ${\mathcal O}(1)$ size solution, in finite time.

In order to prove \eqref{eq:perturbation:lower} we use the following classical idea: the initial data being an eigenfunction of $L$, will grow exponentially in time, while the nonlinear term may be suitably controlled iso that it cannot cancel this  exponential growth. To make this idea concrete, it is convenient to define $L_\delta = L - (\lambda + \delta) I$, where  $0 <  \delta \leq \lambda/8$. Due to the shift, the resolvent of $L_{\delta}$ contains the full right-half of the complex plane, and being a bounded perturbation of the Laplacian we also have that $L_\delta$ generates an analytic semigroup over $L^2$, and $L_\delta^{-1}$ is smoothing by two derivatives.  Then from the Duhamel formula and the interpolation bound $\| L_\delta^{-7/8} N(\theta^\eps) \|_{L^2} \leq C \| \theta^\eps\|_{L^2}^{9/8} \| \Delta \theta^\eps\|_{L^2}^{7/8}$, we obtain the estimate
\begin{align} 
\Vert \theta^{\eps} (t) \Vert_{L^{2}} 
&\leq \| \eps e^{Lt}\phi \|_{L^2} + \int_0^t e^{(\lambda + \delta)(t-s)} \|L_\delta^{7/8} e^{L_\delta(t-s)}\|_{{\mathcal L}(L^2)} \| L_\delta^{-7/8} N(\theta^\eps(s))\|_{L^2} ds \notag \\
&\leq  C_\phi \eps e^{\lambda t} 
+ C \int_0^t e^{(\lambda + \delta)(t-s)} \frac{1}{(t-s)^{7/8}} \| \theta^\eps(s)\|_{L^2}^{9/8} \| \Delta \theta^\eps(s)\|_{L^2}^{7/8} ds.
\label{eq:Duhamel:2}
\end{align}
The key step which allows us to close the estimates is proving the existence of a constant $C^* = C^*(\Th_0,\phi,S)$ such that
\begin{align} \label{eq:higher:regularity}
\| \Delta \theta^\eps \|_{L^2} \leq C^* < \infty
\end{align}
for all sufficiently small $\eps$. The proof of estimate \eqref{eq:higher:regularity} relies on the $C^\alpha$ bound which may be obtained from Theorem~\ref{thm:DeGiorgi} and a suitable bootstrap procedure to reach the $H^2$ regularity level (see \cite[Appendix A]{FriedlanderVicol11b}). Lastly, let $R>C_\phi:= \|\phi\|_{L^2}$ and $T = T(R,\eps)$ be the maximal time such that  
\begin{align*}
\|\theta^\eps(t)\|_{L^2} \leq \eps R e^{\lambda t} \quad \mbox{for all} \quad 0 \leq t \leq T. 
\end{align*}
From   \eqref{eq:Duhamel:2} and \eqref{eq:higher:regularity} one may then obtain
\begin{align} 
\|\theta^\eps(t)\|_{L^2} \leq C_\phi \eps e^{\lambda t} + C_1\left( \eps R e^{\lambda t} \right)^{9/8} \label{eq:theta:eps}
\end{align}
and then in turn that 
\begin{align*}
 T \geq T_\eps = \lambda^{-1} \log(C_2/\eps).
\end{align*}
A lower bound on $\|\theta^\eps(t)\|_{L^2}$, similar to \eqref{eq:theta:eps} but with the linear term minus the $9/8$ power term, may be obtained from the Duhamel formula. In turn, this implies that $\|\theta^\eps(T_\eps)\|_{L^2} \geq C_2(2 C_\phi - R) =: C_* > 0$ and concludes the proof of Theorem~\ref{thm:linear:nonlinear}.

\subsubsection{The existence of unstable eigenvalues for the linearized MG operator}
\label{sec:unstable:eigenvalue:dissipative}
In order to apply Theorem~\ref{thm:linear:nonlinear} to the (nonlinear) MG equations \eqref{eq:MG:1}--\eqref{eq:MG:2}, we need to find a steady state $\Th_0$ for which the corresponding linear operator $L$ defined in \eqref{eq:L}, has unstable spectrum. 

To achieve this we use an idea that dates back at least to~\cite{MeshalkinSinai61} for Navier-Stokes and~\cite{FriedlanderStraussVishik97} for Euler. First, note that $\Th_0 = F(x_3)$, with $F$ smooth, is automatically a steady state, with associated velocity $U_0=0$, if we assume that the source term is given by $S = -\eps_\kappa F''(x_3)$. We thus consider the source term to equal $S = \eps_\kappa m^2 \sin(m x_3)$, for which we have a steady state solution of \eqref{eq:MG:1} given by
\begin{align}
\Th_0 = \sin (m x_3), \quad U_0 = 0. \label{eq:the:steady:state}
\end{align}
For the above choice of the steady states, the equation $\partial_t \theta = L \theta$  becomes
\begin{align}
\partial_t \theta(x,t) = L\theta(x,t) = - m u_3(x,t)  \cos(m x_3) + \eps_\kappa \theta(x,t), \quad u_3 = M_3[\theta], \label{eq:eigenvalue}
\end{align}
where $M_3$ is the Fourier multiplier defined earlier in \eqref{eq:M3}--\eqref{eq:D}.
Next, we make the {\em ansatz} that the temperature $\theta$ which solves \eqref{eq:eigenvalue} has the form
\begin{align}
\theta(x,t) = e^{\sigma t} \sin(k_1 x_1) \sin (k_2 x_2) \sum_{p\geq 1} c_p \sin(m p x_3) \label{eq:ansatz}
\end{align}
for some fixed $k_1, k_2$ in $\ZZ_\ast$, to be chosen precisely later.
Using \eqref{eq:M3}--\eqref{eq:D} we may compute from the third component of the velocity field $u$ associated with $\theta$ given by \eqref{eq:ansatz}, insert this into \eqref{eq:eigenvalue} and obtain that $\sigma$ and $\{c_p\}_{p\geq 1}$ must solve 
\begin{align}
\sigma_p \sum_{p\geq 1} c_p \sin(m p x_3) + \sum_{p\geq 1} \frac{c_p}{\alpha_p} \left( \sin(m(p+1)x_3) + \sin(m (p-1) x_3) \right) = 0 \label{eq:ansatz:eqn}
\end{align}
where we have denoted
\begin{align}
\sigma_p = \sigma + \eps_\kappa \left( k_1^2 + k_2^2 + (mp)^2\right) \label{eq:sigma:p}
\end{align}
and 
\begin{align}
\frac{2}{\alpha_p} =  m N^2 \frac{k_2^2(k_1^2 + k_2^2)}{N^4 (mp)^2 (k_1^2 + k_2^2 + (mp)^2) + k_2^4} \label{eq:alpha:p}
\end{align}
with $N$ being the non-dimensional constant in \eqref{eq:parameters}. An essential feature of the coefficients $\alpha_p$ is that they form an increasing sequence which grows rapidly, as $p^4$ when $p \to \infty$. Equation \eqref{eq:ansatz:eqn} gives a recurrence relation for the sequence $\{c_p\}_{p\geq 1}$ (see e.g. \cite{FriedlanderVicol11b,FriedlanderRusinVicol12}):
\begin{align}
&\sigma_p c_p + \frac{c_{p+1}}{\alpha_{p+1}} + \frac{c_{p-1}}{\alpha_{p-1}} = 0, \quad \mbox{for } p\geq 2, \label{eq:recursion:1}\\
&\sigma_1 c_1 + \frac{c_2}{\alpha_2} = 0, \mbox{for } p=1 \label{eq:recursion:2}.
\end{align}
We emphasize that given any $\sigma > 0$ (which then uniquely defines $\sigma_p$ for all $p \geq 1$) and given any $c_1>0$, one can solve for all $\{c_p\}_{p\geq 2}$ using the recursion relations \eqref{eq:recursion:1}--\eqref{eq:recursion:2}. The key is however that only for suitable values of $\sigma$ do the $c_p$ vanish sufficiently fast so that $\theta$ defined by \eqref{eq:ansatz} is sufficiently smooth. Solving for such $\sigma$ amounts to solving the infinite continued fraction equation
\begin{align}
\sigma_1 \alpha_1 = \frac{1}{\sigma_2 \alpha_2 - \frac{1}{\sigma_3 \alpha_3 - \frac{1}{\sigma_4 \alpha_4 - \ldots}}} \label{eq:continued:frac}.
\end{align}
Solving \eqref{eq:continued:frac} is the main difficulty in this argument, since we need not just show there is a solution, but also give bounds on it in terms of the $\alpha_p$. We omit these details and summarize the results obtained in \cite[Lemma 2.6]{FriedlanderVicol11b}, \cite[Lemma 4.2]{FriedlanderRusinVicol12}, and \cite[Theorem 2.2]{FriedlanderGancedoSunVicol12}  as follows.
\begin{lemma}[\bf Bound of solution to continued fraction equation]\label{lem:contd:frac}
Let $\alpha_p,\sigma_p$ be defined by \eqref{eq:sigma:p} and \eqref{eq:alpha:p} respectively, and let $c_1 = \alpha_1$.  
There exits a real eigenvalue $\sigma_\ast = \sigma_\ast(k_1,k_2,m,N,\eps_\kappa) >0$ which solves \eqref{eq:continued:frac}, and we have the estimate
\begin{align}
\frac{1}{\alpha_1\alpha_2} \leq \sigma_1 \sigma_2 \leq \frac{2}{\alpha_1 \alpha_2}. \label{eq:sigma:est}
\end{align}
Moreover, the corresponding sequence $\{c_p\}_{p\geq 1}$ decays exponentially in $p$.
\end{lemma}
A direct consequence of estimate \eqref{eq:sigma:est} is that the eigenvalue $\sigma_\ast$ obeys the upper and lower bounds
\begin{align}
 & \sigma_\ast \geq \frac{m N^2}{2} \frac{k_2^2(k_1^2 + k_2^2)}{4 N^4 m^2 (k_1^2 + k_2^2 + 4 m^2) + k_2^4}  -\eps_\kappa (k_1^2 + k_2^2 + 4m^2) \label{eq:sigma:est:low}\\
 &\sigma_\ast \leq m N^2  \frac{k_2^2(k_1^2 + k_2^2)}{N^4 m^2 (k_1^2 + k_2^2 + m^2) + k_2^4}  -\eps_\kappa (k_1^2 + k_2^2 + m^2).\label{eq:sigma:est:high}
\end{align}
For $N,\eps_\kappa,m$ fixed, we hence obtain from \eqref{eq:sigma:est:high} that there exists at most finitely many frequency pairs $(k_1,k_2) \in \ZZ^2_\ast$ for which $\sigma_\ast > 0$. Moreover, one may show that in the parameter range of the geodynamo, i.e. $0 < \eps_\kappa \ll 1$,  for suitable $(k_1,k_2)$ the left side of \eqref{eq:sigma:est:low} is positive (see \eqref{eq:max:growth} below), thereby proving: 
\begin{theorem}[\cite{FriedlanderVicol11b}]
Associated to the steady state $\Th_0$ defined in \eqref{eq:the:steady:state}, the linearized MG operator $L$ defined in~\eqref{eq:eigenvalue}, has unstable spectrum. 
\end{theorem}

Lastly, note that one may optimize in $(k_{1},k_{2})$, and compute a lower bound for the largest eigenvalue of $L$ in terms of the physical parameters $N, \eps_\kappa$, and of $m$. More precisely, we have that 
\begin{align} 
\sup_{k_1,k_2 \in \ZZ_\ast} \left(  \frac{m N^2}{2} \frac{k_2^2(k_1^2 + k_2^2)}{4 N^4 m^2 (k_1^2 + k_2^2 + 4 m^2) + k_2^4}  -\eps_\kappa (k_1^2 + k_2^2 + 4m^2) \right) \geq  \frac{1}{2^8 \eps_\kappa}   - 8 m^2 \eps_\kappa \label{eq:max:growth},
\end{align}
with the lower bound being attained e.g. for (of course, we need to take integer part)
\[
k_1 = \frac{1}{16 \eps_\kappa} \quad \mbox{and} \quad k_2 = \frac{N \sqrt{m}}{ \sqrt{8\eps_\kappa}}, \quad \mbox{for } 0 < \eps_\kappa\ll 1.
\]
In particular, we note that as $\eps_\kappa \to 0$ the lower bound given by the right hand side of \eqref{eq:max:growth} blows up (this will be important in Sections~\ref{sec:inviscid} and~\ref{sec:frac}). Moreover, for the above choice of $k_1$ and $k_2$, we have $k_1 k_2^{-2} \to 1/(2 N^2 m)$ as $\eps_\kappa \to 0$. This is in precise accordance to the fact that the absolute value of the Fourier multiplier associated to $M$ grows at the fastest rate when $|k_2| = {\mathcal O}(|k_1|^{1/2})$ as $|k_1| \to \infty$. See \eqref{eq:M:max:growth}.

Combining the results in this subsection with Theorem~\ref{thm:linear:nonlinear}, we have thus constructed an example of solutions to the MG equations which displays strong (exponential) growth in $\theta$ and hence in the velocity field $u$ by \eqref{eq:MG:2} and also in the magnetic field $b(x, t)$  through \eqref{eq:1.6}. The rate of this growth can be shown to be proportional to $\exp(t/ 2^8 \eps_\kappa)$ which is consistent with the dynamo instability scenario.

\section{Non-diffusive MG equation}
\label{sec:inviscid}
As discussed in the introduction, in the fluid core the value of $\eps_\kappa$, given by the quotient of $\kappa$ the molecular diffusivity coefficient, and $L\,V$, the product of the typical length scale and typical velocity in the fluid, is very small. 

Therefore, from the physical point of view it is natural to study the non-diffusive (unforced) MG equation
\begin{align} 
& \partial_t \theta + (u \cdot \nabla) \theta = 0 \label{eq:inv:1}\\
& u = M[\theta], \quad \nabla \cdot u = 0 \label{eq:inv:2}
\end{align}
or the corresponding non-diffusive singular active scalar equation when \eqref{eq:inv:2} is replaced by \eqref{eq:ABS:2}, i.e. $u_j = \partial_i  T_{ij} \theta$.
As we shall see below, the system \eqref{eq:inv:1}--\eqref{eq:inv:2} is {\em ill-posed} in Sobolev spaces, in the sense that there is no operator which is Lipschitz in time at $t=0$. 

\subsection{On the (lack of) well-posedness for the MG equation in Sobolev spaces}
\label{sec:analytic}

Before we discuss the issue of ill-posedness, let us briefly discuss why the {\em evenness} of the operator $M$ breaks the classical proof of local existence in Sobolev spaces for the class of singular active scalar equation \eqref{eq:inv:1}--\eqref{eq:inv:2}. The standard energy-approach to local existence in $H^s$ would entail applying the $(-\Delta)^{s/2}$ to \eqref{eq:inv:1}, and then take $L^2$ inner products with $(-\Delta)^{s/2} \theta$. Morally speaking (this is rigorous when $s/2$ is an integer in view of the Leibniz rule) we have only one term which seem to prevent closing the estimates at the $H^s$ level:
\begin{align}
{\mathcal T} = \int (-\Delta)^{s/2} u \cdot \nabla \theta\; (-\Delta)^{s/2} \theta\; dx = \int M_j (-\Delta)^{s/2}  \theta\; \partial_j \theta \; (-\Delta)^{s/2} \theta\; dx   \label{eq:T:bad}
\end{align}
Indeed, the other term where more than $s$ derivatives lie on a single function is $ \int u \cdot \nabla (-\Delta)^{s/2} \theta (-\Delta)^{s/2} \theta dx$, but this term vanishes upon integrating by parts since $\nabla \cdot u = 0$. The only hope to treat the term ${\mathcal T}$ would be to discover a commutator structure. However, since $M$ is {\em not anti-symmetric}, i.e. {\em even} in Fourier space, we {\em cannot} write ${\mathcal T} = - {\mathcal S}$, where 
\begin{align}
{\mathcal S} = \int M_j \left( (-\Delta)^{s/2}  \theta\; \partial_j \theta \right) (-\Delta)^{s/2} \theta\; dx = {\mathcal T} + \int [M_j, \partial_j \theta] (-\Delta)^{s/2} \theta \; (-\Delta)^{s/2} \theta dx = {\mathcal T}  + {\mathcal C} \label{eq:S:good}
\end{align}
and $[\cdot,\cdot]$ is the usual commutator bracket. Instead for the MG, we have ${\mathcal T} =  {\mathcal S}$, and we {\em cannot} explore the fact that ${\mathcal C}$ may be bounded in terms of $\| \theta \|_{H^s}^3$ (if $s>2$). This is the main reason why we do not seem to be able to close estimates at the Sobolev level. 

On the other hand, it was observed in~\cite{ChaeConstantinCordobaGancedoWu12}, see also \cite{FriedlanderVicol11b}, that if we replace $u = M[\theta]$ (cf.~\eqref{eq:inv:2}) by a certain class of singular constitutive laws $u_j = \partial_i T_{ij}\theta$ (cf.~\eqref{eq:ABS:2}), and the matrix $T_{ij}$ is such that 
\[ \langle T_{ij} f , g \rangle = \langle f, T_{ij} g\rangle \] for any $f,g$ and $i,j$ (here $\langle \cdot, \cdot \rangle$ is the usual $L^2$ inner product), then the corresponding singular active scalar equation \eqref{eq:inv:1}--\eqref{eq:ABS:2} is locally well posed in $H^s$ for $s>2$. The main reason is that under these assumptions the operator $\partial_i T_{ij}$ is anti-symmetric, and we are able to write $2 {\mathcal T} = {\mathcal C}$, using the notation of \eqref{eq:T:bad} and \eqref{eq:S:good}. The commutator term ${\mathcal C}$ would enable one estimates to close at the $H^s$ level.

Lastly, we mention that for analytic initial data, the non-diffusive MG equation is indeed locally well-posed in the class of real-analytic functions, since each term in the equation loses at most one derivative. This was proven in \cite[Theorem~2.1]{FriedlanderVicol11b} by using a Cauchy-Kowalewski-type argument, in the spirit of the Gevrey-class regularity proof of~\cite{FoiasTemam89} (see also~\cite{LevermoreOliver97,OliverTiti00a,KukavicaVicol09}). In particular, this result justifies the short time numerical simulations of the MG equation are justified, as they only keep track of finitely many Fourier modes (hence the associated functions are real analytic).

\subsection{Hadamard ill-posedness for the non-diffusive MG equation} 
\label{sec:ill}

A Cauchy problem for a certain partial differential equation is called  {\em well-posed in $X$ in the sense of Hadamard}, if for any initial data in $X$, the problem has a unique solution in $L^{\infty}(0,T;X)$, with $T$ depending only on the initial data in the $X$-norm, and the solution map $Y \mapsto L^{\infty}(0,T;X)$ satisfies some continuity properties, for a suitable space $Y \subset X$.  If one of these properties fail, the Cauchy problem is called {\em ill-posed}. In this section we discuss the failure of the solution map for the non-diffusive MG equation to be {\em Lipschitz continuous} with respect to perturbations in the initial data (in the topology of a Sobolev space $X$) around a specific steady profile $\Th_{0}$. Such results have been recently obtained for several models in fluid dynamics, see e.g.~\cite{Grenier00,BourgainPavlovic08,Renardy09,CheskidovShvidkoy10,GerardVaretDormy10,GerardVaretNguyen10,GuoNguyen10,DaiSchonbek11,FriedlanderGancedoSunVicol12,FriedlanderRusin12} and references therein.

The proof consists of a linear and a nonlinear step. We construct a steady state and a sequence of eigenfunctions, with arbitrarily large eigenvalues, for the linearized MG equation around this steady state. Once these eigenvalues are exhibited, one may use a fairly robust argument (see, e.g.~\cite{Renardy09}) to show that this severe linear ill-posedness implies the Lipschitz ill-posedness for the nonlinear problem. We emphasize that the key feature of the MG constitutive law \eqref{eq:MG:2} which gives this ill-posednes is the evenness and unboundedness of the symbol associated to $M$. In fact, the evenness also gives non-uniqueness for for bounded weak solutions of the non-diffusive MG equation, as was recently established by Shvydkoy~\cite{Shvydkoy11}.

\subsubsection{Linear ill-posedness}
\label{sec:ill:linear}
Let $m\geq 1$ be a positive integer. We use the steady state 
\begin{align}\label{eq:steady:state}
\Th_0 = \sin(m x_3),
\end{align} 
with associated velocity field $U_0 = 0$. Similarly to \eqref{eq:eigenvalue}, we consider the associated linearized (non-diffusive) MG operator defined by
\begin{align} 
L \theta = - M_3[\theta] \partial_3 \Th_0 \label{eq:steady:L}
\end{align}
where the $M_3$ operator is defined in terms of its Fourier symbol \eqref{eq:M3}--\eqref{eq:D}. With this notation, we have
\begin{theorem}[\cite{FriedlanderVicol11b}] \label{thm:ill:lin}
The Cauchy problem associated to the linearized $MG$ equation 
\begin{align}
\partial_t \theta = L \theta \label{eq:steady:linear}
\end{align}
is ill-posed in the sense of Hadamard over $L^2$. That is, for any $T,K>0$, there exits a real-analytic initial datum $\theta(0,\cdot)$, such that the IVP associated to \eqref{eq:steady:linear} does not have a solution $\theta \in L^\infty(0,T;L^2)$, which obeys
\begin{align} 
\sup_{t\in[0,T]} \| \theta(t,\cdot) \|_{L^2} \leq K \|\theta(0,\cdot) \|_{Y} \label{eq:linear:well}
\end{align}
where $Y$ is any Sobolev space that embeds compactly in $L^2$.
\end{theorem}
The proof of the above theorem builds on the construction of unstable eigenvalues in the critically diffusive case, more specifically on Lemma~\ref{lem:contd:frac}. 

As in Section~\ref{sec:unstable:eigenvalue:dissipative}, cf.~\eqref{eq:ansatz} we make the ansatz that $\theta$ is given in terms of a Fourier sine series 
\[ \theta(t,x) = e^{\sigma t} \sin(k_1 x_1) \sin (k_2 x_2) \sum_{p\geq 1} c_p \sin(m p x_3), \]
with $k_1$ and $k_2$ to be chosen precisely later. Note that since we are in the non-diffusive case, we need to set $\eps_\kappa=0$, and hence cf.~\eqref{eq:sigma:p} we have $\sigma_p = \sigma$ for any $p\geq 1$. The value of $\alpha_p$ remain however unchanged. Using Lemma~\ref{lem:contd:frac}, we obtain an eigenvalue $\sigma_\ast = \sigma_\ast(k_1,k_2)$, and an exponentially decaying sequence $\{c_p\}_{p\geq 1}$, such that cf.~\eqref{eq:sigma:est:low} the bound
\begin{align} 
\sigma_\ast = \sigma_\ast(k_1,k_2) \geq  \frac{m N^2}{2} \frac{k_2^2(k_1^2 + k_2^2)}{4 N^4 m^2 (k_1^2 + k_2^2 + 4 m^2) + k_2^4} \label{eq:sigma:est:non-diffusive:low}
\end{align}
holds.
We set $k_1= j^2$ and $k_2 =j$ in \eqref{eq:sigma:est:non-diffusive:low}, where $j \geq 2m$ is a positive integer,
and obtain the lower bound
\begin{align} 
\sigma_{\ast,j} := \sigma_\ast(j^2,j) \geq j \cdot \frac{m N^2}{24 N^4 m^2  + 2} = \frac{j}{C_1}\label{eq:sigma:est:non-diffusive}
\end{align}
on the eigenvalue, associated to this choice of $(k_1,k_2)$. Here $C_1 = C_1(m,N)$. We also denote by $\theta^{(j)}$ the eigenfunction constructed in this way. Using an additional argument, which requires suitably bounding the exponentially decaying sequence $\{c_p\}_{p\geq 1}$ which defines our eigenfunction, one may in fact additionally normalize $\theta$ and prove (cf.~\cite[Lemma 2.7]{FriedlanderVicol11b}, \cite[Lemma 4.3]{FriedlanderRusinVicol12}) that 
\begin{align} 
\| \theta^{(j)}(0,\cdot) \|_{Y} = 1, \quad \mbox{and} \quad \| \theta^{(j)}(0,\cdot) \|_{L^2} \geq \frac{1}{C_2 j^2},\label{eq:theta:IC:bounds}
\end{align}
where $C_2 = C_2(m,N,Y)$ is a sufficiently large constant. For example, if $Y=H^s$ with $s>0$, we may explicitly compute a suitable such constant.

We conclude the proof of Theorem~\ref{thm:ill:lin} by using a contradiction argument. If such a solution would exist, for a fixed $T$ and $K$, we pick $j = j(T,K,m,N,Y)$ large enough so that
\[
\frac{\exp(T j / (2 C_1))}{C_2 j^2} \geq 2 K
\]
holds. The linear problem~\eqref{eq:steady:linear} can be shown to have a unique $L^2$ solution, and hence by construction we have
\[ \theta^{(j)}(t,x) = \theta^{(j)}(0,x) \exp(\sigma_{\ast,j} t), \]
the assumed estimate \eqref{eq:linear:well} is contradicted by our choice of $j$ and the lower bound \eqref{eq:sigma:est:non-diffusive}. This concludes the proof of Theorem~\ref{thm:ill:lin}.

\subsubsection{Nonlinear ill-posedness}
\label{sec:ill:nonlinear}

Before stating our main ill-posedness result, let us recall  the definition of Lipschitz well-posedness.

\begin{definition}[\cite{GuoNguyen10,FriedlanderVicol11b}] \label{def:well}
Let $Y \subset X \subset W^{1,4}$ be Sobolev spaces. The Cauchy problem for the non-diffusive MG equation 
 is locally Lipschitz $(X,Y)$ {\em well-posed}, if there exist continuous functions $T,K: [0,\infty)^{2} \rightarrow (0,\infty)$,  so that for every pair of initial data $\theta^{(1)}(0,\cdot), \theta^{(2)}(0,\cdot) \in Y$ there exist unique solutions $\theta^{(1)}, \theta^{(2)} \in L^{\infty}(0,T;X)$ of the initial value problem associated to  \eqref{eq:inv:1}--\eqref{eq:inv:2}, that additionally satisfy
\begin{align}
\Vert \theta^{(1)}(t,\cdot) -\theta^{(2)}(t,\cdot) \Vert_{X} \leq K \Vert \theta^{(1)}(0,\cdot) - \theta^{(2)}(0,\cdot) \Vert_{Y}\label{eq:def:well}
\end{align}
for every $t\in [0,T]$, where $T = T(\Vert \theta^{(1)}(0,\cdot) \Vert_{Y}, \Vert \theta^{(2)}(0,\cdot)\Vert_{Y})$ and $K  = K (\Vert \theta^{(1)}(0,\cdot)\Vert_{Y},\Vert \theta^{(2)}(0,\cdot)\Vert_{Y})$.
\end{definition}
We note that Defintion~\ref{def:well} allows the solution map to lose regularity. Since the MG equation is first order, the typical pairs of spaces $(X,Y)$ that we have in mind here are $X = H^{s}$, and  $Y=H^{s+1}$, with $s>1 + 3/4$. 

The main result of this section is then:
\begin{theorem}[\cite{FriedlanderVicol11b}] \label{thm:nonlinear:ill}
The non-diffusive MG equations \eqref{eq:inv:1}--\eqref{eq:inv:2}, are locally Lipschitz $(X,Y)$ {\em ill-posed} in Sobolev spaces $Y\subset X$ embedded in $W^{1,4}$, in the sense of Definition~\ref{def:well}.
\end{theorem}
The proof of the above result is based on the severe linear ill-posedness obtained in Theroem~\ref{thm:ill:lin}, and a standard perturbative argument~\cite[Theorem 2]{Renardy09} or~\cite[pp.~183]{Tao06}, which we breifly sketch below.

Let $\Th_0 = \sin (m x_3)$ be the steady state defined above in \eqref{eq:steady:state}, let $L \theta = - M_3[\theta] \partial_r \Th_0$ be its associated linearized MG operator, and denote the MG nonlinearity by $N(\theta) = M_j [\theta] \partial_j \theta$. By our assumptions on $X$ and $Y$, we have 
\begin{align}
\| L \theta\|_{L^2} \leq C_1 \|\theta\|_X  \quad \mbox{and} \quad \|N(\theta)\|_{L^2} \leq C_1 \|\theta\|_{X}^2 
\label{eq:N:L:bounds}
\end{align}
for some constant $C_1 = C_1(\Th_0) >0$. Also, let $C_0 = \|\Th_0\|_{Y}$.

In definition \eqref{def:well}, we chose $\theta^{(2)}(0) = \Th_0 = \theta^{(2)}(t)$ be the steady state, and let $\theta^{(1)}(0) = \theta^{(1)}_\delta(0) = \Th_0 + \delta \psi_0$, for some smooth function $\psi_0$ to be chosen later. Assume by contradiction that the problem is Lipschitz locally well-posed in $(X,Y)$. One may then show that there exist $T = T(C_0)$ and $K=K(C_0)$ such that 
\[
\sup_{t\in [0,T]} \|\theta^{(1)}_\delta(t) - \Th_0\|_{X} \leq K \delta
\]
for {\em any} $\delta \in (0,\delta_0(C_0)]$. Writing the solution $\theta^{(1)}_\delta(t)$ as an order $\delta$ perturbation around the steady state, $\theta^{(1)}(t) = \Th_0 + \delta \psi_\delta(t)$ we thus have
\begin{align}
 \sup_{t\in[0,T]} \|\psi_\delta(t)\|_X \leq K \label{eq:psi:bound:1}
\end{align}
and $\psi_\delta$ solves the Cauchy problem
\begin{align} 
\partial_t \psi_\delta = L \psi_\delta + \delta N(\psi_\delta), \quad \psi_\delta(0) = \psi_0. \label{eq:psi:eqn}
\end{align}
The key observation is that when combined with \eqref{eq:N:L:bounds}, estimate \eqref{eq:psi:bound:1} yields
\begin{align*} 
\delta \|N(\psi_\delta)\|_{L^2} \leq C_1 \delta K^2
\end{align*}
for any $t \in [0,T]$. Therefore, when passing $\delta \to 0$, the nonlinearity in \eqref{eq:psi:eqn} disappears, and we should obtain the linear limiting equation. Indeed, from \eqref{eq:psi:bound:1}, \eqref{eq:psi:eqn}, and the bounds \eqref{eq:N:L:bounds}, we obtain the necessary estimates to apply the Aubin-Lions compactness argument and show that $\psi_\delta \to \psi$ strongly in $L^\infty(0,T;L^2)$ and on this time interval the limiting function obeys
\begin{align} 
\partial_t \psi = L \psi , \quad \psi(0) = \psi_0. \label{eq:psi:eqn:2}
\end{align}
But $\psi$ also inherits the bound \eqref{eq:psi:bound:1} from the sequence $\psi_\delta$, which contradicts Thereom~\ref{thm:ill:lin} -- choose $\psi_0$ an eigenfunction of $L$ as constructed in Section~\ref{sec:ill:linear}, so that in time $T/2$ the exponential growth of the $L^2$ norm caused the solution to exceed the level $2K$. This concludes our sketch of the proof of Theorem~\ref{thm:nonlinear:ill}.

\section{Fractionally diffusive MG equation and well-prepared initial data}
\label{sec:frac}

The fractionally diffusive MG equation, although non-physical, presents a very challenging and interesting mathematical problem. Replacing the Laplacian in (\ref{eq:MG:1}) by $-(-\Delta)^\gamma$, for $\gamma \in (0,1)$, we obtain the active scalar equation 
\begin{align}
	&\partial_t \theta + u\cdot \nabla \theta = \eps_\kappa (-\Delta)^\gamma \theta + S \label{eq:MGf:1} \\
	& u = M[\theta] \label{eq:MGf:2},\\
	& \theta(0,\cdot) = \theta_0 \label{eq:MGf:3}
\end{align}
where the symbol of the Fourier multiplier $M$ is given by (\ref{eq:M1})--(\ref{eq:D}). Note that here we work with the non-dimensional variables, so that the initial data and the source term are {\em normalized to be of size unity}. 

The main feature of this situation is that the structure of $\widehat{M}_j$ (anisotropy, evenness, and unboundedness) and a relatively weak smoothing effect of the nonlocal operator $(-\Delta)^\gamma$ produce a sharp dichotomy across the value $\gamma = 1/2$, as follows. If $1/2 < \gamma < 1$ the equation is {\em locally} well-posed, while it is Hadamard {\em Lipschitz ill-posed} for $0<\gamma<1/2$. At the threshold value $\gamma = 1/2$ the problem is globally well-posed for $\eps_\kappa \gg 1$, respectively ill-posed for $\eps_\kappa\ll 1$. Furthermore, the anisotropy of the symbols $\widehat{M}_j$ can be used to improve the regularity of the solutions when the initial data and source have a suitably chosen Fourier support. For such {\em well-prepared} data the local existence and uniqueness of solutions can be obtained for all values $\gamma \in (0,1)$, and the global existence holds for all initial data when $\gamma \in (1/2,1)$. We remark that in all above mentioned regimes, the problem is {\em super-critical} in the sense described in Section \ref{sec:diffusive}. 

\subsection{Regime of well-posedness in Sobolev spaces}
\label{sec:frac:well}

In Section \ref{sec:diffusive} we stated the well-posedness result when the diffusion in the equation is given by the full Laplacian. If the diffusion in the equation is introduced by the fractional operator $(-\Delta)^\gamma$, $\gamma \in (1/2,1)$, then we have the following result.

\begin{theorem}[\cite{FriedlanderRusinVicol12}] \label{thm:fractional_local}
	Let $\gamma \in (1/2,1)$, and fix $s > 5/2 + (1-2\gamma)$. Assume that $\theta_0 \in H^s(\TT^3)$ and $S \in L^\infty(0,\infty;H^{s-\gamma}(\TT^3))$ have zero-mean on $\TT^3$. Then there exists a time $T>0$ and a unique smooth solution $ \theta \in L^\infty(0,T; H^s(\TT^3)) \cap L^2(0,T; H^{s + \gamma}(\TT^3))$ of the Cauchy problem (\ref{eq:MGf:1})--(\ref{eq:MGf:2}).
\end{theorem}

Under an additional assumption on the size of $\eps_\kappa$ (equivalently, the size of the initial data in the unscaled equations), we obtain a global existence result.

\begin{theorem}[\cite{FriedlanderRusinVicol12}]\label{thm:fractional_global}
	Let $\gamma$ and $S$ be as in the statement of Theorem~\ref{thm:fractional_local}, and let $ \theta_0 \in H^{s}(\TT^3)$ have zero-mean on $\TT^3$, where $s > 5/2+(1-2\gamma)$. There exists a large enough number  $\bar C$, such that if   $\eps_\kappa \geq \bar C$, then the unique smooth solution $\theta$ of the Cauchy problem \eqref{eq:MGf:1}--\eqref{eq:MGf:2} is global in time, i.e. $\theta \in L^\infty(0,\infty;H^s(\TT^3))$.
\end{theorem}

\begin{remark} \label{rem:rescale}
Note that in~\cite{FriedlanderRusinVicol12} we did not use the non-dimensional variables with which we are working here. Therefore, the global well-posedness condition stated in~\cite{FriedlanderRusinVicol12} is that the initial data and source, have to be {\em sufficiently small}, compared to the diffusivity coefficient $\kappa$. In this paper however, the data and source have been normalized to be of order unity, and the only coefficient remaining is $\eps_\kappa$. It is hence natural that the small data condition turns to a {\em large enough} $\eps_\kappa$ condition. This remains valid throughout this section.
 \end{remark}

The detailed proofs of Theorems~\ref{thm:fractional_local} and ~\ref{thm:fractional_global} are given in \cite{FriedlanderRusinVicol12}, and rely mainly on energy methods and the use of suitable commutator estimates. We omit these details here.

\subsection{Regime of ill-posedness in Sobolev spaces}
\label{sec:frac:ill}

In this subsection we consider the stability of steady states to \eqref{eq:MGf:1}--\eqref{eq:MGf:2}. The construction follows the outline described in Section~\ref{sec:unstable:eigenvalue:dissipative}. In this case, we choose $\Th_0 =  \sin (m x_3)$ and $S = \eps_\kappa m^{2\gamma} \sin(m x_3)$, for some integer $m \geq 1$. 
 As before, the first step in proving that the Cauchy problem \eqref{eq:MGf:1}--\eqref{eq:MGf:2} is Lipschitz ill-posed is to consider the linearized problem around the given steady state, and exhibit the existence of unstable eigenvalues with arbitrarily large real part. Thus, we consider
\begin{align}\label{eq:fractional_inst1}
	\partial_t \theta = L\theta,
\end{align}
where the linear operator $L$ is given by
\begin{align*}
	L\theta(x,t) = -M_3\theta(x,t) \partial_3 \theta_0(x_3) - \eps_\kappa (-\Delta)^{\gamma}\theta.
\end{align*}
We seek the temperature $\theta$ which solves (\ref{eq:fractional_inst1}) and has the form
\begin{align}\label{eq:fractional_theta}
	\theta(t,x) = e^{\sigma t} \sin(k_1 x_{1}) \sin(k_2 x_{2}) \sum_{p\geq 1} c_{p} \sin(m p x_{3})
\end{align}
for some fixed $k_1,k_2 \in \ZZ_*$ to be chosen later. Inserting (\ref{eq:fractional_theta}) into (\ref{eq:fractional_inst1}) and using the definition of $M_3$, we obtain that \eqref{eq:ansatz:eqn} holds, but with $\sigma_{p}$ denoted by
\begin{align}\label{eq:4.8}
	\sigma_{p} = \sigma + \eps_\kappa (k_1^{2} + k_2^{2} + (mp)^{2})^{\gamma}
\end{align}
and with $\alpha_p$ given by \eqref{eq:alpha:p} above. As before, these assumptions lead to solving the continued fraction  equation \eqref{eq:continued:frac} for some real $\sigma >0$. This is achieved by suitably modifying Lemma~\ref{lem:contd:frac} to obtain:
\begin{lemma} \label{lemma:SN:5.2}
 Let $\alpha_{p}, \sigma_p$ be defined by \eqref{eq:alpha:p} and \eqref{eq:4.8} respectively, and let $c_{1} = \alpha_{1}$. There exists a real eigenvalue $\sigma_*>0$, and a rapidly decaying sequence $\{c_{p}\}_{p\geq 2}$ so that the function $\theta$ defined by \eqref{eq:fractional_theta} solves \eqref{eq:fractional_inst1}. Furthermore, we have the lower bound 
 \begin{align} 
\sigma_* \geq  \frac{mN^2}{2} \frac{k_2^2(k_1^2 + k_2^2)}{4N^4 m^2 (k_1^2 + k_2^2 + 4m^2)+ k_2^4} - \eps_\kappa(k_1^2 + k_2^2+ 4m^2)^\gamma. \label{eq:sigma:lower:frac:1}
\end{align}
\end{lemma}

In particular, setting for instance $k_1 = j^2$ and $k_2 = j $ in \eqref{eq:sigma:lower:frac:1}, for each $j \geq 2m$ we obtain a real eigenvalue $\sigma_{\ast,j}$ of $L$, which is bounded from below as
\begin{align}
\sigma_{\ast,j} \geq j^{2\gamma} \left( j^{1-2\gamma} \frac{m N^2}{24 N^4 m^2 + 2} - 3 \eps_\kappa \right).
\label{eq:sigma:lower:frac}
\end{align}
Moreover, given a Sobolev space $Y$ that is compactly embedded in $L^2$, we may ensure the existence of a constant $C_2= C_2(m,N,Y,\gamma)$ such that the smooth eigenfunction $\theta^{(j)}$ corresponding to $\sigma_{\ast,j}$ obeys \eqref{eq:theta:IC:bounds}. The immediate consequence of estimate \eqref{eq:sigma:lower:frac} is that due to $1-2\gamma>0$, we have
\[
\sigma_{\ast,j} \to \infty, \quad \mbox{as} \quad j \to \infty.
\]
The above construction proves that the linearized problem \eqref{eq:fractional_inst1} is ill-posed in $L^2$, as we have exhibited eigenvalues of arbitrarily large real part. Following the argument outlined in Section~\ref{sec:ill:nonlinear}, this type of strong ill-posedness for the linearized problem, automatically implies that the nonlinear problem is Lipschitz ill-posed in Sobolev spaces (recall Defintion~\ref{def:well} of Lipschitz well-posedness). Thus we obtain:

\begin{theorem}[\cite{FriedlanderRusinVicol12}] \label{thm:nonlinear:ill:frac}
The Cauchy problem \eqref{eq:MGf:1}--\eqref{eq:MGf:3}, with $\gamma \in (0,1/2)$ is locally Lipschitz $(X,Y)$ {\em ill-posed} in Sobolev spaces $Y\subset X$ embedded in $W^{1,4}(\TT^3)$. 
\end{theorem}

\subsection{The threshold case \texorpdfstring{$\gamma = 1/2$}{gamma=1/2}}

The value $\gamma = 1/2$ marks the threshold for proving the local well-posedness in Sobolev spaces for \eqref{eq:MGf:1}--\eqref{eq:MGf:3}. If the initial datum is small with respect to $\kappa$, or equivalently, if $\eps_\kappa$ is sufficiently large (see Remark~\ref{rem:rescale} above), then one may use energy estimates to show that there exists a unique global smooth solution evolving from this data. In this direction we have:
\begin{theorem}[\cite{FriedlanderRusinVicol12}] \label{thm:gamma=1/2:small}
Let $s>5/2$ and assume that $\theta_0 \in H^s(\TT^3)$ and $S \in L^\infty(0,\infty;H^{s-\gamma}(\TT^3))$ have zero-mean on $\TT^{3}$. There exists a sufficiently large number $\bar C$, such that  if $\eps_\kappa \geq \bar C$, then  the unique smooth solution $\theta \in L^\infty(0,T;H^s(\TT^3))$ of the Cauchy problem \eqref{eq:MGf:1}--\eqref{eq:MGf:3} is global in time.
\end{theorem}

As in the case of $\gamma \in (1/2,1)$, the proof is based on the classical Kato-Ponce commutator estimate. The proof however does not apply for the case of ``small'' $\eps_\kappa$ (large data), not even to prove the local existence of solutions. 

In the case where $\eps_\kappa$ is suitably small, we may construct a steady state such that solutions evolving from initial data which is arbitrarily close to this steady state diverge from it at an arbitrarily large rate, i.e. the problem is Hadamard ill-posed. More precisely, we have:
\begin{theorem}[\cite{FriedlanderRusinVicol12}] \label{thm:ill:1/2}
 Let $Y \subset X \subset W^{1,4}(\TT^3)$ be Sobolev spaces, let $\gamma=1/2$, and fix the integer $m = 1$. 
Consider the steady state solution of \eqref{eq:MGf:1} given by $\theta_0 ^{(1)}= \sin(m x_3) \in Y$, with source term $S=\eps_\kappa m \sin(m x_3)$, and assume that 
\begin{align}
 \eps_\kappa < \frac{1}{26\sqrt{3}} \label{eq:SN:cond:frac}.
\end{align}
Then, given any $T>0$ and any $K>0$, there exists an initial condition $\theta_0^{(2)} \in Y$ and a corresponding solution $\theta^{(2)} \in L^\infty(0,T;X)$ of the Cauchy problem \eqref{eq:MGf:1}--\eqref{eq:MGf:2}, such that 
 \begin{align*}
\sup_{t\in [0,T]} \Vert  \theta^{(2)}(t,\cdot) - \theta^{(1)}(t,\cdot) \Vert_{X} \geq 2 K \Vert  \theta^{(2)}_0 - \theta^{(1)}_0 \Vert_{Y}
\end{align*}
where $ \theta^{(1)}(t,\cdot) = \theta^{(1)}_0$. That is, the problem \eqref{eq:MGf:1}--\eqref{eq:MGf:2} is locally Lipschitz (X,Y) ill-posed when $\eps_\kappa \ll 1$.
\end{theorem}
To prove the above result, we follow the ill-posedness result in the case $\gamma \in (0,1/2)$. The main difference is that while for $\gamma \in (0,1/2)$ by \eqref{eq:sigma:lower:frac} we have real eigenvalues of arbitrarily large real part, when $\gamma =1/2$ the lower bound \eqref{eq:sigma:lower:frac} only blows up for $\eps_\kappa$ sufficiently small, depending on $N^2$ and our choice of $m$. This condition is rather unsatisfactory, since it implies $\eps_\kappa \to 0$ as $N^2 \to 0$. In fact, we can establish a more robust smallness condition on $\eps_\kappa$, in the sense that it is independent of $N^2$ (and $m$). To see this, insert $k_1 = j^2 N^{-2}$ and $k_2 = j m^{1/2}$ in the lower bound \eqref{eq:sigma:lower:frac:1}, where $j \geq \max\{ N^4 m, 4m\}$. We  obtain the lower bound
\begin{align}
\sigma_{\ast,j} \geq \frac{mN^2}{2} \frac{j^6 m N^{-4}}{12 N^4 m^2  j^4 N^{-4} + j^4 m^2 } - \sqrt{3} \eps_\kappa j^2 N^{-2} = j^2 N^{-2} \left( \frac{1}{26} - \sqrt{3} \eps_\kappa \right) \label{eq:sigma:lower:2}.
\end{align}
Therefore, under the assumption \eqref{eq:SN:cond:frac}, which is independent of $m$ or $N$, the right side of  \eqref{eq:sigma:lower:2} diverges to $\infty$ as $j \to \infty$  and hence we again have a sequence of unstable eigenvalues of the linear problem, with arbitrarily large real part. Combined with the nonlinear argument presented in Section~\ref{sec:ill:nonlinear} this concludes the sketch of the proof of Theorem~\ref{thm:ill:1/2}.

\subsection{Well-prepared initial data} 
\label{sec:frac:prepared}

In this subsection we present results which take advantage of the curved nature (paraboloid-type) of the regions in frequency space where the symbol $\hat{M}(k)$ grows like $|k|$. We observe that if the support of $\hat{\theta}$ lies on a plane, then the operator $M$ behaves like an operator of order zero acting on this $\theta$, and hence the velocity $u$ is as smooth as $\theta$. This allows us to  improve the well-posedness results in the range $\eps_\kappa>0$ and $0 < \gamma < 1$, over the generic case when no conditions on the Fourier transforms of the initial data (and source term) are imposed.

To be more precise let us introduce some notation. For a rational number $q \in \QQ$ with $q \neq 0$, define frequency plane $\mathcal{P}_q$ as
\begin{align*}
	\mathcal{P}_q = \{k = (k_1,k_2,k_3) \in \ZZ^3: k_2 = qk_1 \}.
\end{align*}
One may show that whenever $q \neq 0$, 
\begin{align} \label{eq:M:plane:bound}
|\hat{M}_j(k)| \leq C_{q} < \infty, \quad \mbox{for all } k \in \PP_q \mbox{ and }j \in \{1,2,3\}.
\end{align}
This observation, combined with the fact that when the data $\theta_0$ and source $S$ have frequency support on $\PP_q$, then so does the solution $\theta(t,\cdot)$ of the Cauchy problem \eqref{eq:MGf:1}--\eqref{eq:MGf:3} at all later times $t>0$, means that the smoothing properties of the fractional diffusion term are stronger than they are for the generic data situation. Hence it is possible to prove stronger regularity results, and in particular the {\em local existence} and uniqueness of Sobolev solutions holds {\em for all} values of $\gamma \in (0,1)$.   

\begin{theorem}[\cite{FriedlanderRusinVicol12}] \label{thm:well-prepared:local}
Let $\gamma \in (0,1)$, and fix $s>5/2-\gamma$. Assume that $\theta_0 \in H^s(\TT^3)$ and $S \in L^\infty(0,\infty;H^{s-\gamma}(\TT^3))$ have zero-mean on $\TT^{3}$ and are such that $\supp(\hat\theta_{0}) \subset \PP_{q} \mbox{ and } \supp(\hat{S}(t,\cdot)) \subset  \PP_{q}$ for some fixed $q \in {\mathbb Q}\setminus\{0\}$ and all $t\geq 0$. Then there exists a $T>0$ and a unique smooth solution 
$ \theta \in L^\infty(0,T;H^s(\TT^3)) \cap L^2(0,T;H^{s+\gamma}(\TT^3))$ of the Cauchy problem  \eqref{eq:MGf:1}--\eqref{eq:MGf:3}, such that $\supp(\hat\theta (t,\cdot)) \subset \PP_{q}$ for all $t\in[0,T)$.
\end{theorem}
We note that the above theorem holds even if $\gamma \in (0,1/2)$, a setting in which we know that for {\em generic} initial data the problem is ill-posed. Moreover, one may also prove that if the initial data is small with respect to $\kappa$, or equivalently, $\eps_\kappa \gg 1$ in the non-dimensional setting of this manuscript, the solution obtained in Theorem~\ref{thm:well-prepared:local} is global in time.

The main difficulty in the proof of Theorem~\ref{thm:well-prepared:local} is the construction and uniqueness of solutions, which live on the frequency plane $\PP_q$.
This relies on  constructing of a suitable sequence of approximations supported on the Fourier space on $\PP_q$ and showing its convergence using energy methods. The main obstacle one has to overcome  is designing an iteration scheme which is both suitable for energy estimates, and is such that in each iteration step the frequency support of the approximation lies on $\PP_q$. More explicitly, a scheme in which $\partial_t \theta_{n+1}$ is given in terms of $v \cdot \nabla \theta_n$ preserves the frequency support on $\PP_q$, however is not good for closing the estimates at the level of Sobolev spaces. On the other hand if we define $\partial_t \theta_{n+1}$ in terms of $v \cdot \nabla \theta_{n+1}$, the energy estimates are now clear but it seems difficult to inductively prove that $\hat\theta_{n+1}$ is supported on $\PP_q$. Therefore, we work with a regularized problem, where we introduce an artificial term $-\epsilon \Delta \theta$. The advantage of considering a stronger dissipative problem is that we can {\em construct} smooth solutions which have frequency support lying on $\PP_q$. For details of the construction we refer the interested reader to \cite{FriedlanderRusinVicol12}. Lastly, we pass to the limit $\epsilon \to 0$ and obtain a solution of \eqref{eq:MGf:1}--\eqref{eq:MGf:3} with desired properties. 

If the initial data $\theta_0$ and the source term $S$ have Fourier support in $\PP_q$, we show that for $\gamma \in [1/2,1)$ the local in time solution may be continued for all time.
\begin{theorem}[\cite{FriedlanderRusinVicol12}] \label{thm:well-prepared:global}
Let $s, \theta_0$ and $S$ be as in the statement of Theorem~\ref{thm:well-prepared:local}. If $\gamma \geq 1/2$, the unique smooth solution $\theta$ of the Cauchy problem  \eqref{eq:MGf:1}--\eqref{eq:MGf:3} is global in time.
\end{theorem}
The proof hinges on the observation that the solution constructed in Theorem \ref{thm:well-prepared:local} satisfies $\supp(\hat \theta) \subset \PP_q$. Then, in view of \eqref{eq:M:plane:bound} and the H\"ormander-Mikhlin theorem we have
$
	\|u\|_{W^{s,p}} \leq C \|\theta\|_{W^{s,p}}
$
for all $2 \leq p < \infty, s \geq 0$, and so the velocity $u$ is as smooth as the scalar $\theta$. Hence, the case $\gamma > 1/2$ is ``sub-critical'' for the a priori controlled quantities, and energy estimates in the spirit of  those used to prove the global regularity of the sub-critical SQG equation~\cite{ConstantinWu99,CordobaCordoba04, Wu04}, may be used to prove Theorem \ref{thm:well-prepared:global}. The case $\gamma = 1/2$ in Theorem~\ref{thm:well-prepared:global} is more subtle, and is analogous with the problem of global regularity for the critical SQG equation~\cite{CaffarelliVasseur10, KiselevNazarovVolberg07,KiselevNazarov09,ConstantinVicol12}.
The anisotropic nature of the symbol $M$, makes it more convenient to apply the De Giorgi-inspired proof of \cite{CaffarelliVasseur10} in the context of the fractional MG equation with $\gamma=1/2$. The main ingredient is that if $\theta \in L^\infty$ and $\supp(\hat \theta) \in \PP_q$, then $u \in BMO$. Hence, considering first the linear problem with divergence-free $BMO$ drift we obtain H\"older continuity of the solution. The later norm is sub-critical and one may bootstrap to higher regularity. We refer the interested reader to~\cite{FriedlanderRusinVicol12} for  details.

\section{Conclusions}
\label{sec:concl}

We have proven that the magneto-geostrophic model proposed by Moffatt and Loper is capable of sustaining dynamo action. The MG-dynamo instability exhibited here relies crucially on the fact that the dimensionless parameter $N^2$~\eqref{eq:parameters} is neither $0$, nor $\infty$. The exponential growth is of order ${\mathcal O}(\exp(t/\eps_\kappa))$, for all $N^2 \in (0,\infty)$, and $0< \eps_\kappa \ll 1$.
In terms of the physics in the model, this condition on $N^2$ requires that {\em all three forces}, Coriolis, Lorentz, and gravity, {\em are non-zero}, making the MG-instability somewhat different in nature from the usual Rayleigh-Taylor instability.

At the level of mathematical analysis several problems remain to be addressed, including some well/ill-posedness issues. For the non-diffusive equation, {\em is there a ``stronger'' form of ill-posedness such as norm inflation or even finite time blowup?} Can we refine the well-posedness results for the fractionally diffusive problem in the range $1/2<\gamma<1$? That is, {\em can one establish the global existence of smooth solutions or on the contrary, prove finite time blowup?}

Several, interesting, mathematical issues arise when returning to the derivation of the MG equation from the full three dimensional MHD equations in the spirit of~\cite{MoffattLoper94,Moffatt08}. In particular, we note that Moffatt's and Loper's initial discussion includes the kinematic viscosity $\nu$, which is directly proportional to $\eps$ in \eqref{eq:1.1}. The values of $\eps$ and $\eps_\kappa$ are highly conjectural, but both are believed to be {\em extremely} small. 
Therefore, as {\em opposed} to the MG system, if we keep $\eps>0$ we derive from \eqref{eq:1.1}--\eqref{eq:1.4} the system
\begin{align}
& N^2   {e}_3 \times u = - \nabla P +  {e}_2 \cdot \nabla b + N^2 \theta  {e}_3 + \eps \Delta u \label{eq:C:1.1} \\
& 0 =  {e}_2 \cdot \nabla u + \Delta b \label{eq:C:1.2} \\
& \partial_t \theta + u \cdot \nabla \theta = \eps_\kappa \Delta \theta + S \label{eq:C:1.3} \\
& \nabla \cdot u = 0, \nabla \cdot b = 0. \label{eq:C:1.4}
\end{align}
The presence of the Laplacian term in \eqref{eq:C:1.1} makes the constitutive law relating $u$ and $\theta$ ``better'', in the sense that the $M$ operator corresponding to \eqref{eq:C:1.1}--\eqref{eq:C:1.4}, now depends both on $N$ and $\eps$, and is in fact smoothing of degree $2$ (see equations (29)--(30) in \cite{MoffattLoper94}). Hence, in this case the constitutive law is even better than that for the two-dimensional Euler equation in vorticity form, implying that: {\em in the range $\eps>0$, the system \eqref{eq:C:1.1}--\eqref{eq:C:1.4} is globally well-posed, even if $\eps_\kappa = 0$}. 

This observation is interesting in view of the much more challenging problem of {\em establishing properties of the solutions to \eqref{eq:C:1.1}--\eqref{eq:C:1.4} in the vanishing parameters limits $\eps \to 0$ and $\eps_\kappa \to 0$}, as proposed by Moffatt and Loper. The later question is directly related to the turbulent scaling regime for $\theta$ conjectured by Moffatt in~\cite{Moffatt08}. ``Turbulence'' in the framework of vanishing diffusivity $\eps_\kappa$ for the forced system \eqref{eq:C:1.1}--\eqref{eq:C:1.4}, and the presence of anomalous dissipation of energy for $\theta$, or ``thetargy'', may be plausible for very weak solutions (below Onsager critical) of the non-diffusive equation.

\subsection*{Acknowledgements} The work of SF is supported in part by the NSF grants DMS-0803268 and DMS-1207780. The work of VV is supported in part by the NSF grant DMS-1211828.

\newcommand{\etalchar}[1]{$^{#1}$}


\begin{thebibliography}{CCC{\etalchar{+}}12}

\bibitem[BGS02]{BardosGuoStrauss02}
C.~Bardos, Y.~Guo, and W.~Strauss.
\newblock Stable and unstable ideal plane flows.
\newblock {\em Chinese Ann. Math. Ser. B}, 23(2):149--164, 2002.
\newblock Dedicated to the memory of Jacques-Louis Lions.

\bibitem[BP08]{BourgainPavlovic08}
J.~Bourgain and N.~Pavlovi{{\'c}}.
\newblock Ill-posedness of the {N}avier-{S}tokes equations in a critical space
  in 3{D}.
\newblock {\em J. Funct. Anal.}, 255(9):2233--2247, 2008.

\bibitem[BR95]{BraginskyRoberts95}
S.I. Braginsky and P.H. Roberts.
\newblock Equations governing convection in earth's core and the geodynamo.
\newblock {\em Geophysical \& Astrophysical Fluid Dynamics}, 79(1-4):1--97,
  1995.

\bibitem[CC04]{CordobaCordoba04}
A.~C{\'o}rdoba and D.~C{\'o}rdoba.
\newblock A maximum principle applied to quasi-geostrophic equations.
\newblock {\em Comm. Math. Phys.}, 249(3):511--528, 2004.

\bibitem[CCC{\etalchar{+}}12]{ChaeConstantinCordobaGancedoWu12}
D.~Chae, P.~Constantin, D.~C{{\'o}}rdoba, F.~Gancedo, and J.~Wu.
\newblock Generalized surface quasi-geostrophic equations with singular
  velocities.
\newblock {\em Communications on Pure and Applied Mathematics},
  65(8):1037---1066, 2012.

\bibitem[CCW01]{ConstantinCordobaWu01}
P.~Constantin, D.~Cordoba, and J.~Wu.
\newblock On the critical dissipative quasi-geostrophic equation.
\newblock {\em Indiana Univ. Math. J.}, 50(Special Issue):97--107, 2001.
\newblock Dedicated to Professors Ciprian Foias and Roger Temam (Bloomington,
  IN, 2000).

\bibitem[CMT94]{ConstantinMajdaTabak94}
P.~Constantin, A.J. Majda, and E.~Tabak.
\newblock Formation of strong fronts in the {$2$}-{D} quasigeostrophic thermal
  active scalar.
\newblock {\em Nonlinearity}, 7(6):1495--1533, 1994.

\bibitem[CS10]{CheskidovShvidkoy10}
A.~Cheskidov and R.~Shvydkoy.
\newblock Ill-posedness of the basic equations of fluid dynamics in besov
  spaces.
\newblock {\em Proc. Amer. Math. Soc}, 138(3):1059--1067, 2010.

\bibitem[CV10]{CaffarelliVasseur10}
L.A. Caffarelli and A.~Vasseur.
\newblock Drift diffusion equations with fractional diffusion and the
  quasi-geostrophic equation.
\newblock {\em Ann. of Math. (2)}, 171(3):1903--1930, 2010.

\bibitem[CV12]{ConstantinVicol12}
P.~Constantin and V.~Vicol.
\newblock Nonlinear maximum principles for dissipative linear nonlocal
  operators and applications.
\newblock {\em Geometric And Functional Analysis}, 22(5):1289--1321, 2012.

\bibitem[CW99]{ConstantinWu99}
P.~Constantin and J.~Wu.
\newblock Behavior of solutions of 2{D} quasi-geostrophic equations.
\newblock {\em SIAM J. Math. Anal.}, 30(5):937--948, 1999.

\bibitem[CW08]{ConstantinWu08}
P.~Constantin and J.~Wu.
\newblock Regularity of {H}{\"o}lder continuous solutions of the supercritical
  quasi-geostrophic equation.
\newblock {\em Ann. Inst. H. Poincar{\'e} Anal. Non Lin{\'e}aire},
  25(6):1103--1110, 2008.

\bibitem[DG57]{DeGiorgi57}
E.~De~Giorgi.
\newblock Sulla differenziabilit{\`a} e l'{a}naliticit{\`a} delle estremali
  degli integrali multipli regolari.
\newblock {\em Mem. Accad. Sci. Torino. Cl. Sci. Fis. Mat. Nat. (3)}, 3:25--43,
  1957.

\bibitem[DQS11]{DaiSchonbek11}
M.~Dai, J.~Qing, and M.E. Schonbek.
\newblock Norm inflation for incompressible magneto-hydrodynamic system.
\newblock {\em arXiv preprint arXiv:1110.2723}, 2011.

\bibitem[FGSV12]{FriedlanderGancedoSunVicol12}
S.~Friedlander, F.~Gancedo, W.~Sun, and V.~Vicol.
\newblock On a singular incompressible porous media equation.
\newblock {\em Journal of Mathematical Physics}, 53(11):1--20, 2012.

\bibitem[Fil12]{Filonov12}
N.~Filonov.
\newblock On the regularity of solutions to the equation $- {\Delta} u+ b
  \nabla u= 0$.
\newblock {\em arXiv preprint arXiv:1209.1243}, 2012.

\bibitem[FPS06]{FriedlanderPavlovicShvydkoy06}
S.~Friedlander, N.~Pavlovi{{\'c}}, and R.~Shvydkoy.
\newblock Nonlinear instability for the {N}avier-{S}tokes equations.
\newblock {\em Comm. Math. Phys.}, 264(2):335--347, 2006.

\bibitem[FPV09]{FriedlanderPavlovicVicol09}
S.~Friedlander, N.~Pavlovi{{\'c}}, and V.~Vicol.
\newblock Nonlinear instability for the critically dissipative
  quasi-geostrophic equation.
\newblock {\em Comm. Math. Phys.}, 292(3):797--810, 2009.

\bibitem[FR12]{FriedlanderRusin12}
S.~Friedlander and W.~Rusin.
\newblock On the second iterate for active scalar equations.
\newblock {\em arXiv preprint arXiv:1206.3349}, 06 2012.

\bibitem[FRV12]{FriedlanderRusinVicol12}
S.~Friedlander, W.~Rusin, and V.~Vicol.
\newblock On the supercritically diffusive magneto-geostrophic equations.
\newblock {\em Nonlinearity}, 25(11):3071--3097, 2012.

\bibitem[FSV97]{FriedlanderStraussVishik97}
S.~Friedlander, W.~Strauss, and M.~Vishik.
\newblock Nonlinear instability in an ideal fluid.
\newblock {\em Ann. Inst. H. Poincar{\'e} Anal. Non Lin{\'e}aire},
  14(2):187--209, 1997.

\bibitem[FT89]{FoiasTemam89}
C.~Foias and R.~Temam.
\newblock Gevrey class regularity for the solutions of the {N}avier-{S}tokes
  equations.
\newblock {\em J. Funct. Anal.}, 87(2):359--369, 1989.

\bibitem[FV11a]{FriedlanderVicol11a}
S.~Friedlander and V.~Vicol.
\newblock Global well-posedness for an advection-diffusion equation arising in
  magneto-geostrophic dynamics.
\newblock {\em Ann. Inst. H. Poincar\'e Anal. Non Lin\'eaire}, 28(2):283--301,
  2011.

\bibitem[FV11b]{FriedlanderVicol11b}
S.~Friedlander and V.~Vicol.
\newblock On the ill/well-posedness and nonlinear instability of the
  magneto-geostrophic equations.
\newblock {\em Nonlinearity}, 24(11):3019--3042, 2011.

\bibitem[FV12]{FriedlanderVicol12}
S.~Friedlander and V.~Vicol.
\newblock Higher regularity of {H}{\"o}lder continuous solutions of parabolic
  equations with singular drift velocities.
\newblock {\em J. Math. Fluid Mech.}, 14(2):255--266, 2012.

\bibitem[GC87]{GhilChildress87}
M.~Ghil and S.~Childress.
\newblock {\em Topics in geophysical fluid dynamics: atmospheric dynamics,
  dynamo theory, and climate dynamics}.
\newblock Springer-Verlag New York Inc., New York, NY, 1987.

\bibitem[GN10]{GuoNguyen10}
Y.~Guo and T.~Nguyen.
\newblock A note on the prandtl boundary layers.
\newblock {\em arXiv:1011.0130v3 [math.AP]}, 2010.

\bibitem[GOC04]{GlatzmaierOgdenClune04}
G.A. Glatzmaier, D.E. Ogden, and T.L. Clune.
\newblock Modeling the earth's dynamo.
\newblock {\em Geophysical Monograph Series}, 150:13--24, 2004.

\bibitem[GR97]{GlatzmaierRoberts97}
G.A. Glatzmaier and P.H. Roberts.
\newblock Simulating the geodynamo.
\newblock {\em Contemporary physics}, 38(4):269--288, 1997.

\bibitem[Gre00]{Grenier00}
E.~Grenier.
\newblock On the nonlinear instability of {E}uler and {P}randtl equations.
\newblock {\em Comm. Pure Appl. Math.}, 53(9):1067--1091, 2000.

\bibitem[GVD10]{GerardVaretDormy10}
D.~G\'{e}rard-Varet and E.~Dormy.
\newblock On the ill-posedness of the {P}randtl equation.
\newblock {\em J. Amer. Math. Soc.}, 23(2):591--609, 2010.

\bibitem[GVN10]{GerardVaretNguyen10}
D.~G\'{e}rard-Varet and T.~Nguyen.
\newblock Remarks on the ill-posedness of the prandtl equation.
\newblock {\em arXiv:1008.0532v1 [math.AP]}, 2010.

\bibitem[IKR12]{IgnatovaKukavicaRyzhik12}
M.~Ignatova, I.~Kukavica, and L.~Ryzhik.
\newblock The harnack inequality for second-order elliptic equations with
  divergence-free drifts.
\newblock {\em Preprint}, 2012.

\bibitem[Kis10]{Kiselev10}
A.~Kiselev.
\newblock Regularity and blow up for active scalars.
\newblock {\em Math. Model. Nat. Phenom.}, 5(4):225--255, 2010.

\bibitem[KN09]{KiselevNazarov09}
A.~Kiselev and F.~Nazarov.
\newblock A variation on a theme of {C}affarelli and {V}asseur.
\newblock {\em Zap. Nauchn. Sem. S.-Peterburg. Otdel. Mat. Inst. Steklov.
  (POMI)}, 370(Kraevye Zadachi Matematicheskoi Fiziki i Smezhnye Voprosy Teorii
  Funktsii. 40):58--72, 220, 2009.

\bibitem[KNV07]{KiselevNazarovVolberg07}
A.~Kiselev, F.~Nazarov, and A.~Volberg.
\newblock Global well-posedness for the critical 2{D} dissipative
  quasi-geostrophic equation.
\newblock {\em Invent. Math.}, 167(3):445--453, 2007.

\bibitem[KV09]{KukavicaVicol09}
I.~Kukavica and V.~Vicol.
\newblock On the radius of analyticity of solutions to the three-dimensional
  {E}uler equations.
\newblock {\em Proc. Amer. Math. Soc.}, 137(2):669--677, 2009.

\bibitem[LO97]{LevermoreOliver97}
C.D. Levermore and M.~Oliver.
\newblock Analyticity of solutions for a generalized {E}uler equation.
\newblock {\em J. Differential Equations}, 133(2):321--339, 1997.

\bibitem[ML94]{MoffattLoper94}
H.K. Moffatt and D.E. Loper.
\newblock The magnetostrophic rise of a buoyant parcel in the earth's core.
\newblock {\em Geophysical Journal International}, 117(2):394--402, 1994.

\bibitem[Mof08]{Moffatt08}
H.K. Moffatt.
\newblock Magnetostrophic turbulence and the geodynamo.
\newblock In Y.~Kaneda, editor, {\em I{UTAM} {S}ymposium on {C}omputational
  {P}hysics and {N}ew {P}erspectives in {T}urbulence, Nagoya, Japan, September,
  11--14, 2006}, volume~4 of {\em IUTAM Bookser.}, pages 339--346. Springer,
  Dordrecht, 2008.

\bibitem[MS61]{MeshalkinSinai61}
L.D. Me{\v{s}}alkin and J.G. Sina{\u\i}.
\newblock Investigation of the stability of a stationary solution of a system
  of equations for the plane movement of an incompressible viscous liquid.
\newblock {\em J. Appl. Math. Mech.}, 25:1700--1705, 1961.

\bibitem[MX12]{MiaoXue12}
C.~Miao and L.~Xue.
\newblock Global well-posedness for a modified critical dissipative
  quasi-geostrophic equation.
\newblock {\em Journal of Differential Equations}, 252(1):792--818, 2012.

\bibitem[Nas58]{Nash58}
J.~Nash.
\newblock Continuity of solutions of parabolic and elliptic equations.
\newblock {\em Amer. J. Math.}, 80:931--954, 1958.

\bibitem[NU12]{NazarovUraltseva12}
A.I. Nazarov and N.N. Ural'tseva.
\newblock The harnack inequality and related properties for solutions of
  elliptic and parabolic equations with divergence-free lower-order
  coefficients.
\newblock {\em St. Petersburg Math. J}, 23:93--115, 2012.

\bibitem[Osa87]{Osada87}
H.~Osada.
\newblock Diffusion processes with generators of generalized divergence form.
\newblock {\em J. Math. Kyoto Univ.}, 27(4):597--619, 1987.

\bibitem[OT00]{OliverTiti00a}
M.~Oliver and E.S. Titi.
\newblock Remark on the rate of decay of higher order derivatives for solutions
  to the {N}avier-{S}tokes equations in {${\bf R}^n$}.
\newblock {\em J. Funct. Anal.}, 172(1):1--18, 2000.

\bibitem[Ren09]{Renardy09}
M.~Renardy.
\newblock Ill-posedness of the hydrostatic {E}uler and {N}avier-{S}tokes
  equations.
\newblock {\em Arch. Ration. Mech. Anal.}, 194(3):877--886, 2009.

\bibitem[RG00]{RobertsGlatzmaier00}
P.H. Roberts and G.A. Glatzmaier.
\newblock Geodynamo theory and simulations.
\newblock {\em Reviews of modern physics}, 72(4):1081, 2000.

\bibitem[RG01]{RobertsGlatzmaier01}
P.H. Roberts and G.A. Glatzmaier.
\newblock The geodynamo, past, present and future.
\newblock {\em Geophysical \&amp; Astrophysical Fluid Dynamics},
  94(1-2):47--84, 2001.

\bibitem[Shv11]{Shvydkoy11}
R.~Shvydkoy.
\newblock Convex integration for a class of active scalar equations.
\newblock {\em J. Amer. Math. Soc.}, 24(4):1159--1174, 2011.

\bibitem[SS{\v{S}}Z12]{SereginSilvestreSverakZlatos12}
G.~Seregin, L.~Silvestre, V.~{\v{S}}ver{\'a}k, and A.~Zlato{\v{s}}.
\newblock On divergence-free drifts.
\newblock {\em Journal of Differential Equations}, 252(1):505--540, 2012.

\bibitem[SV12]{SilvestreVicol12}
L.~Silvestre and V.~Vicol.
\newblock H{\"o}lder continuity for a drift-diffusion equation with pressure.
\newblock {\em Annales de l'Institut Henri Poincare (C) Non Linear Analysis},
  20(4):637--652, 2012.

\bibitem[Tao06]{Tao06}
T.~Tao.
\newblock {\em Nonlinear dispersive equations}, volume 106 of {\em CBMS
  Regional Conference Series in Mathematics}.
\newblock Published for the Conference Board of the Mathematical Sciences,
  Washington, DC, 2006.
\newblock Local and global analysis.

\bibitem[VF03]{VishikFriedlander03}
M.~Vishik and S.~Friedlander.
\newblock Nonlinear instability in two dimensional ideal fluids: the case of a
  dominant eigenvalue.
\newblock {\em Comm. Math. Phys.}, 243(2):261--273, 2003.

\bibitem[Wu05]{Wu04}
J.~Wu.
\newblock Global solutions of the 2{D} dissipative quasi-geostrophic equation
  in {B}esov spaces.
\newblock {\em SIAM J. Math. Anal.}, 36(3):1014--1030 (electronic), 2004/05.

\bibitem[Zha04]{Zhang04}
Q.S. Zhang.
\newblock A strong regularity result for parabolic equations.
\newblock {\em Comm. Math. Phys.}, 244(2):245--260, 2004.

\end{thebibliography}
\end{document}